\documentclass[12pt,letterpaper]{article}
\usepackage[top=2cm, bottom=4.5cm, left=2.5cm, right=2.5cm]{geometry}
\usepackage{amsmath,amsthm,amsfonts,amssymb,amscd}
\usepackage{lastpage}
\usepackage{enumerate}
\usepackage[shortlabels]{enumitem}
\usepackage{fancyhdr}
\usepackage{mathrsfs}
\usepackage{tikz}
\usetikzlibrary{arrows}
\usepackage{xcolor}
\usepackage{graphicx}
\usepackage{soul}
\usepackage{listings}
\usepackage{hyperref}
\usepackage{esvect}
\usepackage{stmaryrd}
\usepackage{graphicx}
\usepackage{tcolorbox}
\usepackage{algorithm} 
\usepackage{algpseudocode}
\usepackage{float}
\usepackage{lmodern}
\usepackage[normalem]{ulem}
\usepackage[toc,page]{appendix}
\usepackage{xcolor} 
\usepackage{multirow}
\usepackage{algorithm} 
\usepackage{algpseudocode}
\usepackage{varwidth}
\usepackage{hyperref}
\hypersetup{%
  colorlinks=true,
  linkcolor=black,
  linkbordercolor={0 0 1}
}
\usepackage{braket}
\newcommand{\R}{\mathbb{R}}

\usepackage{accents}

\usepackage{authblk}

\title{A hybrid SIAC -- data-driven post-processing filter for discontinuities in solutions to numerical PDEs}

\author[$\dagger$]{Soraya Terrab}
\author[$\dagger$,$\ddagger$]{Samy Wu Fung}
\affil[$\dagger$]{\footnotesize{Department of Applied Mathematics and Statistics, Colorado School of Mines, USA}}
\affil[$\ddagger$]{\footnotesize{Department of Computer Science, Colorado School of Mines, USA}}
\author[$\ast$]{Jennifer K. Ryan}
\affil[$\ast$]{\footnotesize{Department of Mathematics, KTH Royal Institute of Technology, Sweden}}
\date{\footnotesize{\today}}                    

\begin{document}

\maketitle

\begin{abstract}
     \footnotesize{
We present a hybrid filter that is only applied to the approximation at the final time and allows for reducing errors away from a shock as well as near a shock. It is designed for discontinuous Galerkin approximations to PDEs and combines a rigorous moment-based Smoothness-Increasing Accuracy-Conserving (SIAC) filter with a data-driven CNN filter. While SIAC improves accuracy in smooth regions, it fails to reduce the $\mathcal{O}(1)$  errors near discontinuities, particularly in inviscid compressible flows with shocks. Our hybrid SIAC–CNN filter, trained exclusively on top-hat functions, enforces consistency constraints globally and higher-order moment conditions in smooth regions, reducing both $\ell_2$ and $\ell_\infty$ errors near discontinuities and preserving theoretical accuracy in smooth regions. We demonstrate its effectiveness on the Euler equations for the Lax, Sod, and Shu-Osher shock-tube problems.}
\end{abstract}

\section{Introduction}
\label{sec:intro}
Traditional filtering for numerical approximations is designed to accelerate the decay of the Fourier coefficients away from discontinuities.  However, they are not aimed to improve the approximation near a discontinuity.  In this article, we propose a hybrid filter that effectively handles both smooth and discontinuous regions in numerical solutions of PDEs. This filter only needs to be applied at the final time of the approximation. Many physical applications—ranging from aerospace and automotive engineering to combustion systems—rely on accurate solutions of the Euler equations~\cite{BatinaEulerAircraft, WoodgateHawkAircraft, ElmqvistAutomotive, FrezzottiCombustion, SchulzeGasTurbineCombustors}. However, the presence of shock waves in inviscid flows introduces discontinuities that conventional high-order shock-capturing methods (e.g., discontinuous Galerkin (DG), essentially non-oscillatory (ENO), and weighted essentially non-oscillatory (WENO) schemes) struggle to resolve, often resulting in $\mathcal{O}(1)$ errors near discontinuities.

Post-processing filters, such as the Smoothness-Increasing Accuracy-Conserving (SIAC) filter~\cite{SIACnonuniform, Docampo2017, Picklo2022}, are commonly used to recover accuracy away from discontinuities by post-processing the numerical solution. Such filtering is particularly important when numerical data is affected by aliasing errors, discretization-induced pointwise errors, the Gibbs phenomenon, or other artifacts~\cite{XuliaCAF, HESTHAVEN2016441, Spectral, Picklo2022, Ryan2005, Vandeven1991FamilyOS}. Filters are designed to restore uniform convergence in the approximation; however, their effectiveness is generally limited to smooth regions. If a discontinuity arises in an otherwise $\mathcal{C}^\infty$ smooth solution to a PDE, the piecewise $\mathcal{C}^\infty$ solution cannot be recovered with the same order of accuracy due to global pollution effects introduced by the discontinuity. Mock and Lax~\cite{MockLax} demonstrated that post-processing the numerical approximation allows recovery of the exact solution and its derivatives with an order $\nu-\delta$, where $\delta$ depends on the number of vanishing moments of the post-processing kernel filter, if suitable pre-processing of the data is considered. Subsequent advancements in post-processing kernels have achieved higher accuracy orders, such as order $2\nu -1$ for the SIAC kernel~\cite{Ryan2005} and a user-specified order for the exponential kernel, which does not strictly satisfy the definition of a filter as outlined in~\cite{Vandeven1991FamilyOS}. Despite these advancements, the enhanced accuracy remains confined to smooth regions, while the error in the vicinity of a discontinuity remains at order one~\cite{Spectral}. In regions around a shock or contact discontinuity, the numerical solution can only converge to the average of the left and right limits of the discontinuity~\cite{Spectral, MockLax}.

Several strategies have been developed to address shock resolution in PDE solvers. Traditional shock-capturing methods naturally allow discontinuities to emerge in the solution \cite{LeVeque92, Toro}, yet they incur spurious oscillations due to the Gibbs phenomenon \cite{Spectral}. Post-processing techniques, such as ENO stencil corrections \cite{LafonOsher1990,LafonOsher1992} and the SIAC filter \cite{Ryan2005}, have been employed to recover higher accuracy in smooth regions, though their efficacy near discontinuities remains limited. More recently, machine learning approaches, including physics-informed neural networks (PINNs) \cite{RaissiPINN, MaoPINNEuler, JagtapCPINNs, LiuDiscPINNs, LiangContDiscPINNs}, have shown promise in various PDE contexts but still face challenges in accurately capturing shock behavior.

To overcome the filtering limitations and be able to resolve shocks, our work integrates a data-driven convolutional neural network (CNN) filter — trained exclusively on top-hat discontinuities — with the SIAC filter to post-process discontinuities. This hybrid SIAC–CNN approach enforces consistency constraints about discontinuities and the SIAC moment conditions away from discontinuities, thereby reducing both $\ell_2$ and $\ell_\infty$ errors in discontinuous areas while preserving theoretical accuracy in smooth regions. Convolutional filters based on CNNs have introduced a data-driven approach to post-processing and filtering numerical data. CNNs are widely used in fields such as visual tracking and recognition~\cite{Danelljan2016, Han2022, REHMAN2019171} and have been combined with other neural network (NN) techniques to extract spatial and temporal frequency information from multi-frequency data~\cite{Hu_CNN-LSTM}. While NN kernels have been trained to approximate basis functions for the Gaussian kernel in mean-field control~\cite{Vidal}, and CNNs have been leveraged to learn time-dependent PDEs from limited data~\cite{QU_CNN-PDEs}, the use of CNNs to learn a post-processing kernel for numerical PDE solutions under physical constraints has not been explored.

Our main contribution is a post-processing filter that only needs to be applied at the final time and integrates an output-constrained CNN with the SIAC filter, augmented by Hermite polynomial interpolation in the vicinity of discontinuities. This hybrid SIAC–CNN filter is specifically designed to reduce the $\mathcal{O}(1)$ error near shocks while maintaining the high-order accuracy of SIAC in smooth regions. We demonstrate the effectiveness of our method on benchmark shock-tube problems for the Euler equations (Lax, Sod, and Shu-Osher). The remainder of this paper is organized as follows: we first present  the relevant background (Section~\ref{sec:background}), followed by the development of the hybrid filter (Section~\ref{sec:hybrid-filter}), the experimental setup (Section~\ref{sec:experiments}), and finally the results and discussion (Sections~\ref{sec:results} and~\ref{sec:discussion}).

\section{Background}\label{sec:background}
\subsection{Filters}\label{sec:filters}
We briefly review two key filtering concepts: spectral filtering and the SIAC kernel, which improve numerical approximations away from discontinuities.

\subsubsection{Spectral Filters}\label{sec:spectral}
Consider approximating a PDE solution $u(x,t)$ on $\Omega\subset\R$ with a truncated Fourier series using $M$ modes:
\[
\mathcal{P}_M u(x,t) = \sum_{m=-M/2}^{M/2} \hat{u}_m(t) e^{i\frac{2\pi m}{|\Omega|}x},\quad \text{with }\quad \hat{u}_m(t) = \frac{1}{|\Omega|}\int_\Omega u(x,t)e^{-i\frac{2\pi m}{|\Omega|}x}\,dx.
\]
\par If we consider an approximation with a discontinuity, the error at the discontinuity will remain large, and the error between the exact solution $u$ and the approximation $\mathcal{P}_Mu$ is up to order $1/\sqrt{M}$, even if the resolution is increased. While the function may be smooth away from the discontinuity and periodic, the global rate of convergence is governed by the presence of the discontinuity due to the global expansion basis. The slow convergence away from a discontinuity and the non-uniform convergence near the discontinuity constitute the behavior in Gibbs phenomenon \cite{Spectral}, with Gibbs oscillations polluting the numerical approximation over the entire domain.

\par Filtering is utilized to recover spectral accuracy \textit{away from the discontinuity} \cite{Vandeven1991FamilyOS}. Filtering modifies the Fourier coefficients by a real, even function $\sigma(\eta)$ that satisfies:
\begin{subequations}
\begin{align}
    \sigma(\eta) &= 0 \quad \text{for } |\eta|\ge1,\\[1mm]
    \sigma(0)&=1,\\[1mm]
    \sigma^{(i)}(0)=\sigma^{(i)}(1)&=0,\quad i=1,\dots,q-1.
\end{align}
\end{subequations}
The filtered approximation {is obtained by simply applying the filter to the individual Fourier coefficients:
\[
\mathcal{P}_M^\sigma u(x,t)=\sum_{m=-M/2}^{M/2}\sigma\left(\frac{|m|}{M/2}\right)\hat{u}_m(t) e^{i\frac{2\pi m}{|\Omega|}x}
\]
This allows for faster decay of high-frequency modes and improved convergence in smooth regions \cite{Vandeven1991FamilyOS, Spectral}.

\subsubsection{The SIAC Kernel}\label{sec:siac}
 The SIAC filter was developed as a post-processing tool to extract accuracy and increase smoothness in DG methods for hyperbolic equations \cite{Ryan2005, SIACnonuniform}. This filter is based on the earlier work done by Bramble and Schatz \cite{BrambleSchatz} and Cockburn et al. \cite{CockburnShu2003}. The SIAC filter achieves superconvergence by reducing errors and increasing convergence rates given polynomial reproduction in the SIAC kernel construction. A DG numerical approximation of order $p+1$ can achieve up to $2p+1\ (2\nu-1 = 2p+1)$ order of accuracy in the $L^2$ norm after SIAC post-processing for linear hyperbolic equations~\cite{Ryan2005, XuliaCAF}.  The SIAC filter has been adapted for nonuniform meshes \cite{SIACnonuniform}, multi-dimensional problems using tensor products \cite{SIACnonuniform},  line filtering \cite{Docampo2017, XuliaCAF}, and hexagonal splines \cite{HexSIAC}, as well as enhanced multiresolution analysis \cite{Picklo2022}.

 \par The SIAC kernel is constructed as a weighted average of central B-splines:
\[
K^{(r+1,k+1)}(x) = \sum_{\gamma=0}^{r} c_{\gamma} B^{(k+1)}(x-x_\gamma),
\]
where the coefficients $c_\gamma$ are chosen to enforce consistency ($\int_\R K(x)dx=1$) and polynomial reproduction up to degree $r$. The SIAC filter is applied as a convolution to a numerical approximation $u_h$ with discretization $h= \frac{\Omega}{N}$ where $N$ is the number of elements in the mesh,
\[
u_h^*(x, T) = \frac{1}{H}\int_\R K^{(r+1,k+1)}\left(\frac{x-\xi}{H}\right)u_h(\xi, T)d\xi.
\]
Here, the kernel scaling is $H$, typically equal to the element size, $h$.

\subsection{Euler Equations}\label{sec:Euler}
Our study focuses on the Euler equations for compressible gas dynamics, which are known to develop shock waves. In conservative form, they are written as
\[
\boldsymbol{u}_t + \boldsymbol{f}(\boldsymbol{u})_x = 0,\quad \boldsymbol{u}=(\rho, \rho u, E)^T,\quad \boldsymbol{f}(\boldsymbol{u})=(\rho u, \rho u^2+p, (E+p)u)^T,
\]
where the conservative variables are density $\rho$, momentum $\rho u$, and energy $E$, while the primitive variables are density $\rho$, velocity $u$, and pressure $p$. The corresponding flux terms describe the conservation of mass, momentum, and energy, respectively, in the medium. The total energy is given by
\[
E=\frac{p}{\gamma-1}+\frac{1}{2}\rho u^2,
\]
with the ratio of specific heats at constant pressure and volume is given as $\gamma=1.4$ for an ideal gas. These equations are solved using a modal DG method \cite{Cockburn1999, CockburnShu2001}.

\subsubsection{Initial Conditions}\label{sec:ic}
In the numerical experiments, we consider three standard shock-tube problems in one dimension: Lax~\cite{Lax1954}, Sod~\cite{SOD19781}, and Shu-Osher~\cite{ShuOsher1989}. Each problem is defined on $\Omega=[-5,5]$ with two distinct initial states that generate rarefaction waves, contact discontinuities, and shock waves. While ``shock" is often interpreted as any type of discontinuity, we clarify the distinction between the different discontinuities in these invsicid compressible flow problems, as derived by Jeffrey \cite{Jeffrey}. The contact discontinuity occurs when the jump about the density profile is non-zero, whereas those of the pressure and the fluid velocity vectors are zero, while a shock discontinuity results in discontinuities about the density, velocity, and pressure profiles. Additionally, the rarefaction wave is another discontinuity type. We focus herein on the contact and shock discontinuities that produce the discontinuity in the function itself conducive to the $\mathcal{O}(1)$ error profile that is most challenging to tackle in post-processing methods.

For instance, the Sod problem~\cite{SOD19781} is specified by
\[
(\rho_0, u_0, p_0) = \begin{cases} (1,0,1), & x<0, \\ (0.125,0,0.1), & x\ge 0, \end{cases}
\]
and is evolved until $T_f=2$. The Lax and Shu-Osher problems feature sharper discontinuities, presenting greater challenges due to spurious oscillations. Reference solutions (Fig.~\ref{fig:exact_ref}) computed using the Clawpack solver, as used in \cite{LeVequeRiemann}, illustrate that while rarefaction waves are continuous with discontinuities in the first derivative, shock and contact discontinuities are discontinous, producing the $\mathcal{O}(1)$ error behavior that our filtering approach aims to reduce.

\begin{figure}[t]
    \centering
    \includegraphics[width=0.8\linewidth]{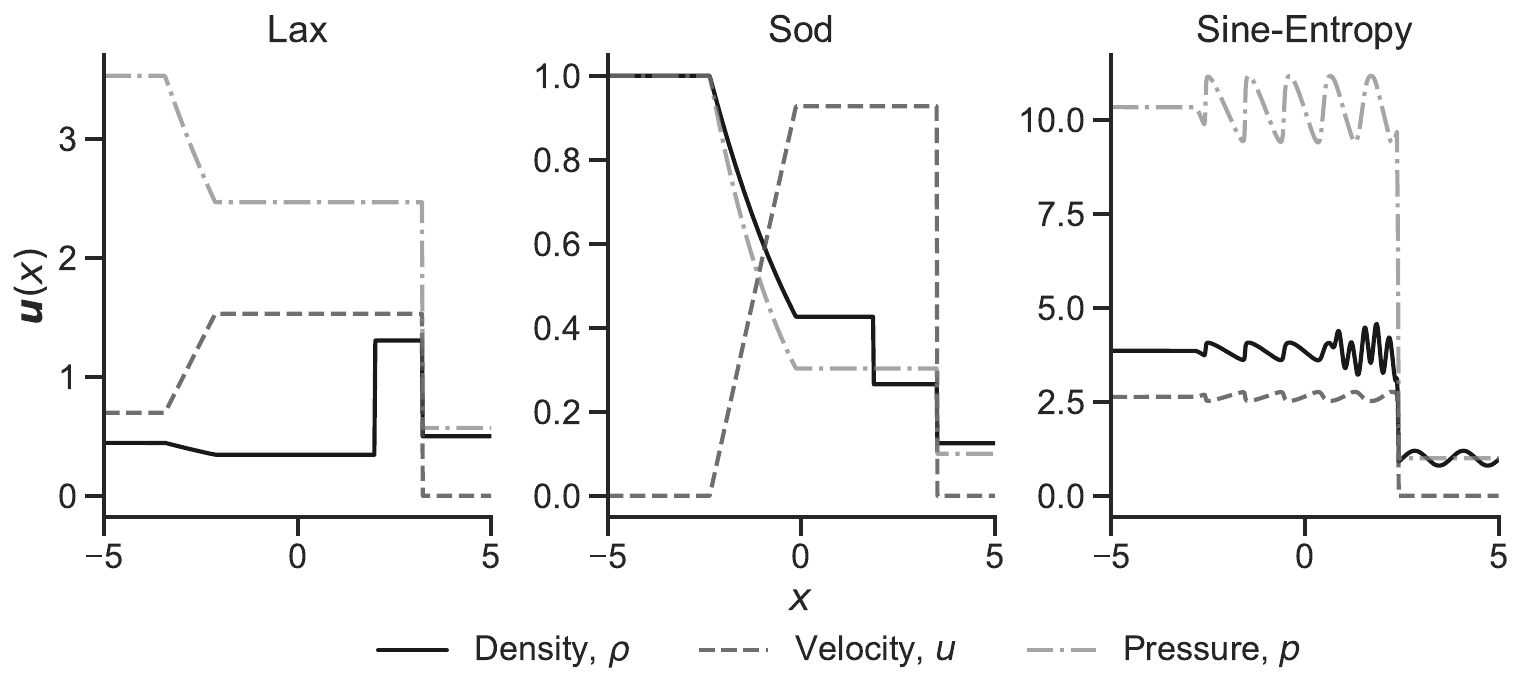}
    \caption{\small{The Euler \protect\hypertarget{refsoln}{exact reference solutions} for the Lax, Sod, and sine-entropy (Shu-Osher) test problems at their respective benchmark final times, $T_f=\{1.3, 2, 1.8\}$. Note the different data ranges for each of the problems. }}
    \label{fig:exact_ref}
\end{figure}

\section{A Hybrid SIAC -- Data-Driven Filter}
\label{sec:hybrid-filter}
\par 
\hypertarget{hybridSIAC} We propose a hybrid, data-driven post-processing filter that is only applied at the final time.  It reduces errors near discontinuities while preserving accuracy in smooth regions. Our approach combines a consistency-only SIAC filter with a data-driven CNN filter and uses interpolation to fuse together their outputs. 

\begin{algorithm}[t]
	\caption{CNN Forward Step} 
	\begin{algorithmic}[1]
        \Function{forward}{$u_h^{\ast}$}: 
            \State compute NN-filtered data: $u_h^{\ast_{NN}} = K_\Theta \ast u_h^{\ast}$
            \State compute normalizing constant: $c_\Theta = \frac{1}{N}\sum K_\Theta \ast \mathbf{1}^N$
            \State return $u_h^{\ast_{NN}}/c_\Theta$
        \EndFunction
	\end{algorithmic} 
    \label{algorithm}
\end{algorithm}

\begin{figure}[t]
    \centering
    \includegraphics[width=0.9\linewidth]{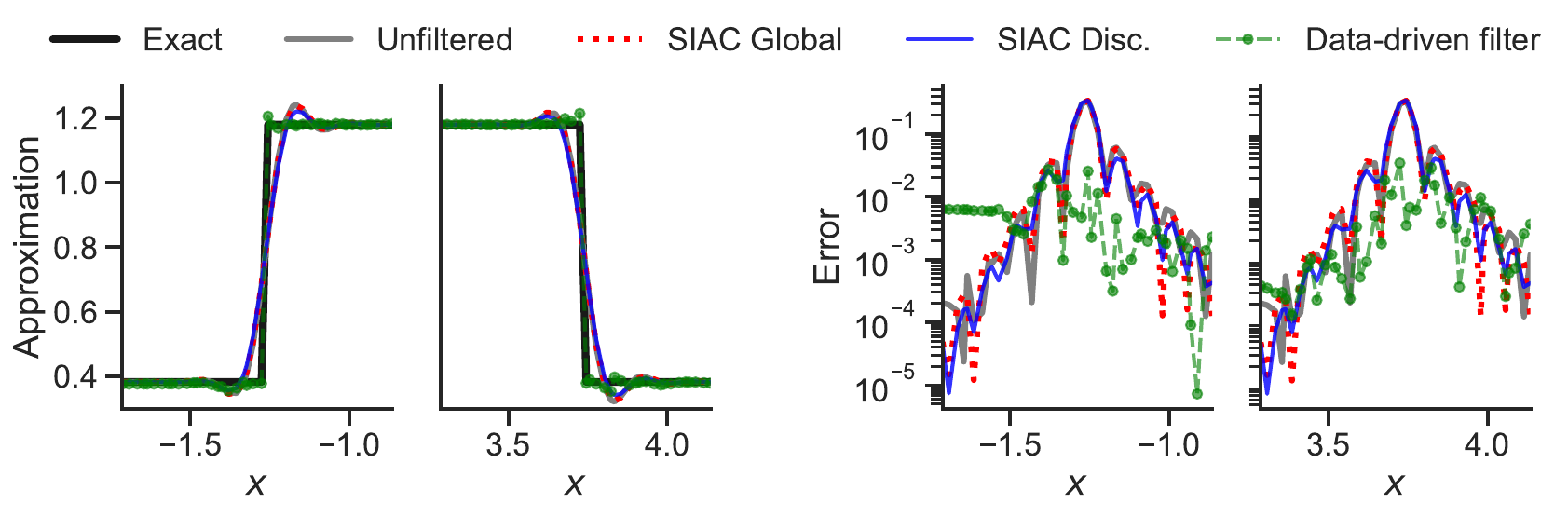}
    \caption{\small{Comparison of filtering approaches: global SIAC (red dotted) with full kernel support and order-2 B-splines, adapted SIAC for discontinuities (blue) with a single order-1 B-spline, and the data-driven filter (green dashed). The left plots show the approximated contact and shock discontinuities for a DG degree-3 sample, while the right plots display the corresponding pointwise errors.}}
    \label{fig:siac_vs_datadriven}
\end{figure}

\subsection{The Data-Driven Filter}\label{sec:data-driven}
Our data-driven filter employs a residual CNN that satisfies a kernel consistency constraint. Specifically, we train a kernel $K_\Theta$ via the constrained optimization problem
\[
\min_\Theta \frac{1}{M} \sum_{j=1}^M \| K_\Theta \ast (K_H \ast u_{j,h}) - u_j \|_2^2 \quad \text{subject to} \quad \int_{\R} K_\Theta(x) \, dx = 1,
\]
where $u_j$ is the exact solution for sample $j$, $u_{j,h}$ is the corresponding DG approximation, and $K_H$ is the SIAC kernel (with scaling $H$) adapted for discontinuity regions to preserve consistency alone. The constraint is enforced by computing the normalizing constant
\[
c_\Theta = \int_{\mathbb{R}} K_\Theta(x)\, dx,
\]
and using automatic differentiation in PyTorch~\cite{paszke2019pytorch} to backpropagate through the normalized output. Algorithm~\ref{algorithm} summarizes the forward step of the CNN filter. More complex constraints can be included using implicit networks~\cite{fung2022jfb,mckenzie2023faster,mckenzie2024three,heaton2021feasibility}, but we leave this for future work.

\subsection{Constructing the Hybrid Filter}\label{sec:construct-hybrid}
The SIAC filter is designed for reducing errors in smooth regions, not discontinuities, as illustrated in Fig.~\ref{fig:siac_vs_datadriven}. Even when using a SIAC filter with reduced kernel support (a moving average kernel), the $\mathcal{O}(1)$ error remains near discontinuities. Our hybrid filter targets the cell with the maximum error by applying the data-driven CNN filter in this region, while Hermite polynomial interpolation fuses its output with the consistency-only SIAC filter applied over the discontinuity window. This approach mitigates over-smoothing from higher-order SIAC filters and preserves pointwise continuity and derivative information in adjacent cells.

\begin{figure}[t]
    \centering
    \includegraphics[width=0.4\linewidth]{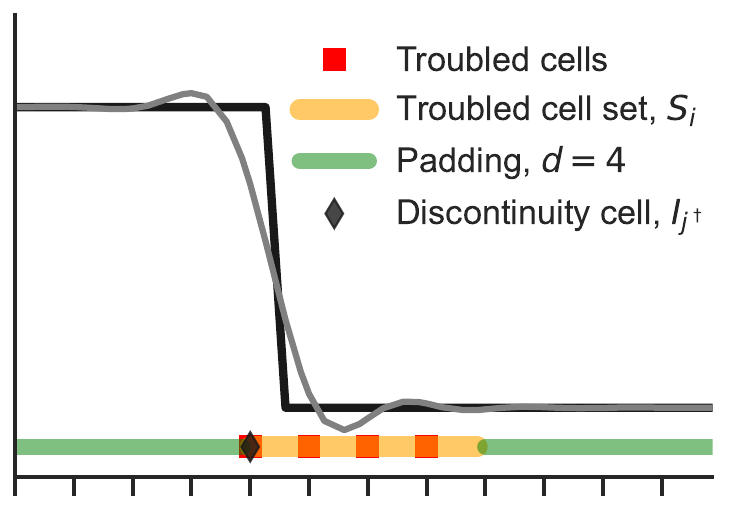}
    \caption{\small{Demonstration of the discontinuity window about a shock discontinuity with an unfiltered DG approximation (gray) and its exact solution (black), where the following are shown: troubled cells (red squares at the left boundary of the cell/element), the grouping of all troubled cells into set $S_i$ (yellow), and the padding of $d=4$ elements (green) to the left and right of set $S_i$, and the cell where the discontinuity is located (black diamond). }}
    \label{fig:disc_window}
\end{figure} 

\subsection{Discontinuity Windows}\label{sec:disc_window}
Accurate identification of discontinuity locations is crucial for our hybrid filter. We use troubled cell indicators to detect discontinuity regions, grouping adjacent troubled cells into sets $S_i$. These sets are padded by $d$ elements on each side to form the discontinuity window $[\min(S_i)-d, \max(S_i)+d]$, as shown in Fig.~\ref{fig:disc_window}. Within this window, the exact cell containing the discontinuity, denoted $I_{j^\dagger}$ (marked by a black diamond), is located using forward differences. The consistency-only SIAC filter is applied over the entire window, while the data-driven filter is applied to $I_{j^\dagger}$ and Hermite interpolation is used for the adjacent cells. This combination minimizes the spread of the $\mathcal{O}(1)$ error and yields a hybrid filtered approximation that reduces global $\ell_\infty$ errors while maintaining high accuracy in smooth regions.

\section{Procedure}\label{sec:experiments}
\par 
\hypertarget{proc}{In this section,} we describe the procedure for building and implementing the hybrid filter, including details on training the data-driven filter (Section~\ref{sec:training}) and constructing the hybrid filtered solution (Section~\ref{sec:hybrid-app}). The code repository is available at \href{https://github.com/sterrab/hybrid-filter}{https://github.com/sterrab/hybrid-filter}.

\subsection{Training the Data-Driven Filter}\label{sec:training}

\begin{figure}[t]
\centering
\includegraphics[width=6.5in]{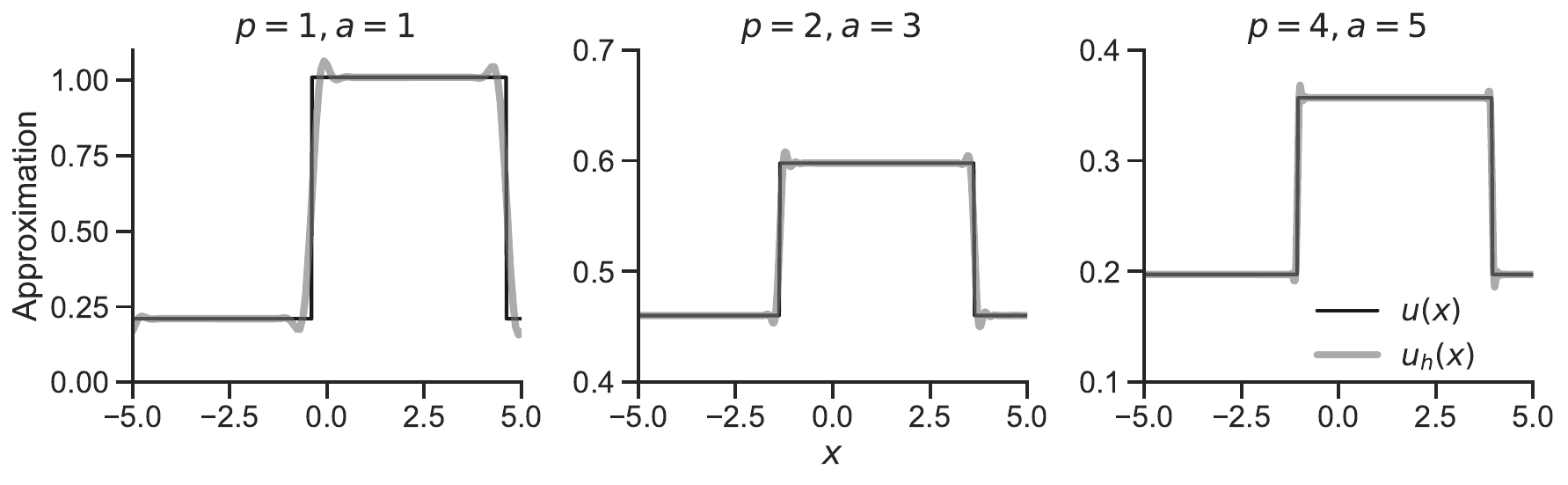}
\caption{\small{Tophat samples with varying DG degree $p$ and wave speed $a$ generated for jump discontinuity data (DG approximations $u_h$ versus their exact solutions $u$).}}
\label{fig:training_samples}
\end{figure}

We train the data-driven filter on synthetic linear data (the linear advection equation $u_t + au_x=0$) with the tophat initial condition
\[
u_0(x) = 
\begin{cases}
\alpha+\delta, & x\in[-2.5,2.5],\\[1mm]
\alpha, & \text{else,}
\end{cases} \quad x\in[-5,5],
\]
using wave speeds $a\in[1,5]$, taken to encompass the range in eigenvalues exhibited by the Sod, Lax, and Shu-Osher shock-tube test problems. We restrict our training set to data from linear equations because of the linearization often involved in approximating solutions to nonlinear equations. This data, generated with a third-order strong stability-preserving Runge-Kutta scheme on a fixed discretization of $N=128$ elements, mimics the discontinuities found in Euler shock waves. Parameters are drawn uniformly from: $\alpha\in[0.1, 0.5)$, $\delta\in[0.1,1)$, DG degree $p\in\{1,2,3,4\}$, and final time $T_f\in \frac{|\Omega|}{a}\cdot[1.1,1.3)$. Examples are shown in Fig.~\ref{fig:training_samples}.

Discontinuity data is extracted from the SIAC-filtered solution $u_h^*=K_H\ast u_h$ using a multiwavelet troubled cell indicator \cite{Vuik2014}. For each flagged troubled cell $TC_i$, a window spanning $[TC_i-4,TC_i+4]$ is defined, as illustrated in the left panel of Fig.~\ref{fig:discontinuity_data}. Data is evaluated at 4 Gauss-Legendre quadrature nodes per element (yielding a 36-dimensional input for the CNN) and normalized by the input data. The right panel of Fig.~\ref{fig:discontinuity_data} shows sorted, normalized samples. In total, the training set comprises 900 top-hat samples (with 100 samples per wave speed) yielding 1,138 discontinuity window samples.

A fixed CNN architecture is used: one input layer, five hidden convolutional layers (kernel size 7, 128 channels, leaky ReLU with slope 0.1), and one output layer. The network is trained with the Adam optimizer (learning rate $10^{-4}$, batch size 200) to minimize the mean-squared error between the CNN-filtered SIAC data and the exact solution, subject to the constraint
\[
\int_{\R} K_\Theta(x)\,dx=1.
\]
Algorithm~\ref{algorithm} outlines the CNN forward step. The best parameters (lowest MSE on a validation set of 50 Euler samples) are saved for $K_\Theta$.

\begin{figure}[t]
    \centering
    \includegraphics[width=0.7\linewidth]{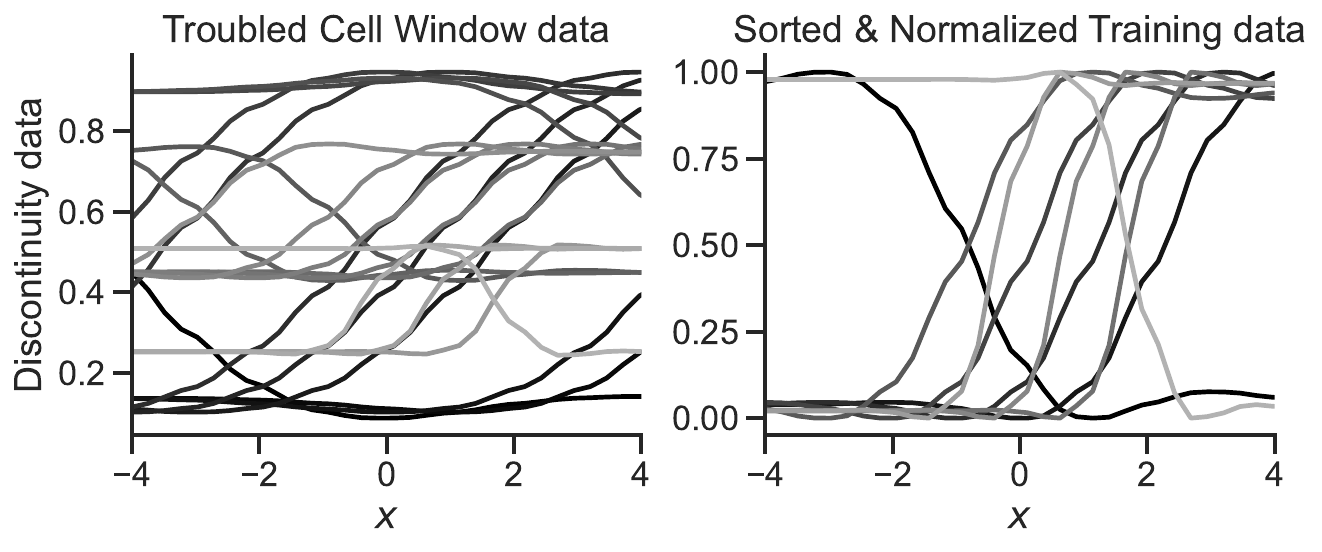}
    \caption{\small{(Left) Window about troubled cell data for cells $[TC_i-4,TC_i+4]$ flagged from samples in Fig.~\ref{fig:training_samples}; (right) Sorted and normalized discontinuity data for NN training.}}
    \label{fig:discontinuity_data}
\end{figure}

\subsection{The Hybrid Filtered Approximation}\label{sec:hybrid-app}

The hybrid filter is applied at the final time and within discontinuity windows identified as groups $S_i$ (cells at most $n=4$ apart) padded by $d=4$ cells, as shown in Fig.~\ref{fig:disc_window}. Within each window, a consistency-only SIAC filter is applied, except for the cell $I_{j^\dagger}$ where the discontinuity is located; there, the data-driven filter is used:
\[
K_\Theta \ast u_h^*\biggr|_{I_{j^\dagger}} = K_\Theta\ast (K_H \ast u_h)\biggr|_{I_{j^\dagger}},
\]
and its two adjacent cells $I_{j^\dagger\pm 1}$ are updated via Hermite polynomial interpolation of degree $p^H=\min(p,2)$ using one point from $I_{j^\dagger}$ and $p^H$ points from $I_{j^\dagger\pm2}$.

\subsubsection{Hybrid-Filtering Euler Data}\label{sec:hybrid-euler}

\begin{figure}[t]
    \centering
    \includegraphics[width=0.95\linewidth]{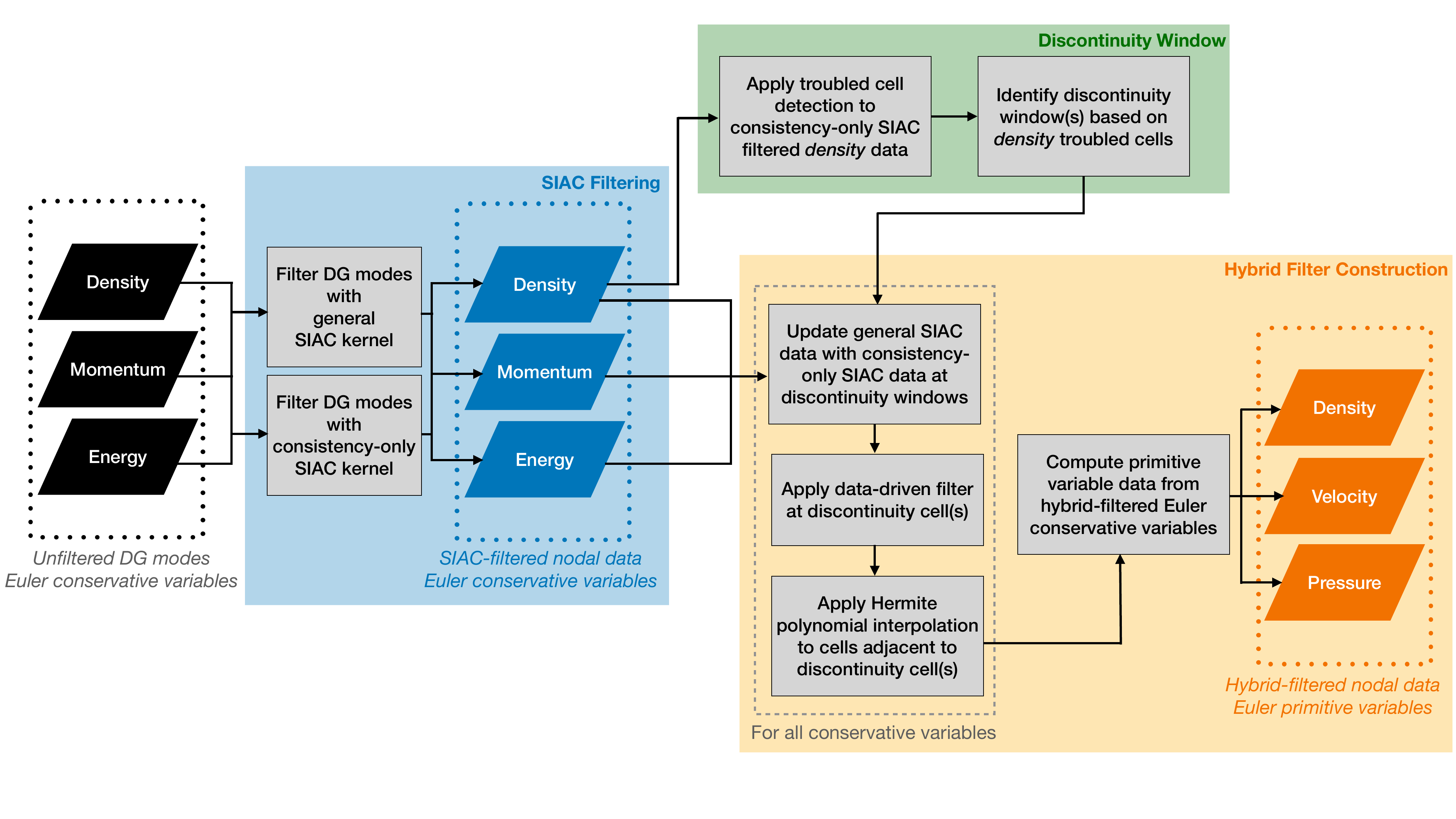}
    \caption{\small{Flowchart of the hybrid filter post-processing of Euler solutions, from unfiltered conservative variables to hybrid-filtered primitive variable outputs.}}
    \label{fig:process}
\end{figure}

For Euler equations, the hybrid filter is applied to the conservative variables, and the filtered primitive variables (density, velocity, pressure) are recovered afterwards. As depicted by the flowchart in Fig.~\ref{fig:process}, SIAC filtering is first applied globally. Then, discontinuity windows are determined using the SIAC-filtered density variable (sufficient for troubled cell indication \cite{Vuik2014}), to which the hybrid filter as presented in Section~\ref{sec:hybrid-app} is applied. The hybrid filtered conservative variables can then be converted to primitive ones to ensure that contact and shock discontinuities are accurately captured while reducing global $\ell_\infty$ errors.

\section{Numerical Results}
\label{sec:results}
\begin{table}[t]
    \caption{\small{The \textbf{Lax} grid $\ell_2$ and $\ell_\infty$ error quartiles for the unfiltered, SIAC moving average, and hybrid filter approximations about discontinuity windows from 84 Lax simulations with final time $T_f\in[1, 1.3]$, resulting in 142 discontinuity window samples $[\min(S_i) - 4,\ \max(S_i) + 4]$. Hybrid filtered error values lower than their unfiltered counterparts appear in bold. \emph{Note:} velocity and pressure profiles are filtered in smooth regions where the density contact discontinuity is flagged as discontinuous.}}
    \centering
    \vspace{3pt}
    \begin{tabular}{|l|l|c|c|c|c|c|c|}
    \hline
        \multirow{2}{*}{Variable} & \multirow{2}{*}{Quartile} & \multicolumn{3}{c|}{$\ell_2$ Error}& \multicolumn{3}{c|}{$\ell_\infty$ Error}\\
        \cline{3-8} 
         & & Unfiltered & SIAC & Hybrid & Unfiltered & SIAC & Hybrid \\
         \hline
         Density & 75\% & 2.01e-01 & 2.12e-01 & \textbf{1.81e-01}  &  4.56e-01 & 4.50e-01 & 4.78e-01 \\
         & Median & 1.68e-01 & 1.74e-01 & \textbf{4.80e-02}  &  4.20e-01 & 4.16e-01 & \textbf{9.51e-02} \\
         & 25\% & 1.48e-01 & 1.57e-01 & \textbf{2.36e-02}  &  3.74e-01 & 3.75e-01 & \textbf{4.96e-02} \\
         \hline 
        Velocity & 75\% & 3.25e-01 & 3.69e-01 & \textbf{8.76e-02}  &  9.01e-01 & 9.70e-01 & \textbf{2.10e-01}\\
         & Median & 2.73e-01 & 3.16e-01 & \textbf{4.65e-02}  &  7.14e-01 & 8.11e-01 & \textbf{9.54e-02} \\
         & 25\% & 9.06e-03 & 4.45e-03 & \textbf{5.19e-03}  &  8.92e-03 & 5.09e-03 & \textbf{6.44e-03}\\
          \hline 
        Pressure & 75\% & 3.41e-01 & 3.68e-01 & \textbf{8.01e-02}  &  9.23e-01 & 9.14e-01 & \textbf{1.39e-01}\\
        & Median & 2.90e-01 & 3.32e-01 & \textbf{3.43e-02}  &  7.52e-01 & 7.64e-01 & \textbf{6.90e-02}\\
        & 25\% & 1.35e-02 & 5.37e-03 & \textbf{6.24e-03}  &  1.57e-02 & 5.89e-03 & \textbf{7.53e-03} \\
        \hline
    \end{tabular}
    \label{tab:error_dataset}
\end{table}

\begin{figure}[t]
    \centering
    \includegraphics[width=\linewidth]{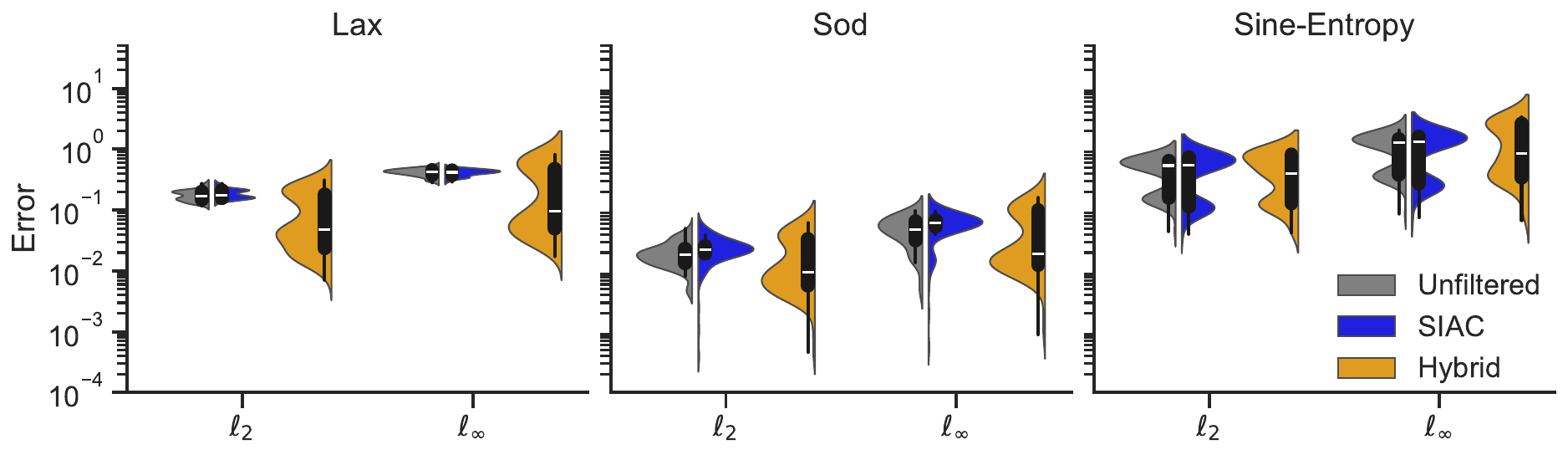}
    \caption{\small{\textbf{Euler density violinplots} for the grid $\ell_2$ and $\ell_\infty$ error in the \textit{(left to right) Lax, Sod, and sine-entropy (Shu-Osher)} datasets. Each violinplot includes a box-and-whisker plot (median in white) and a kernel density estimate for the error distribution.}}
    \label{fig:dataset_density}
\end{figure}

\par \hypertarget{results}{The hybrid filter} is tested on the jump discontinuities of the Lax \cite{Lax1954}, Sod \cite{SOD19781}, and sine-entropy (Shu-Osher) \cite{ShuOsher1989} problems, introduced in Section~\ref{sec:ic}. All simulations (with $N=128$ elements in $\Omega=[-5,5]$) are run to or near benchmark final times. A TVB troubled-cell detector \cite{QiuShu} and Krivodonova moment limiter \cite{KRIVODONOVA2007} keep solutions bounded yet allow some oscillations, thus emulating noisy data for filter testing.

\par The Euler simulations include 84 Lax samples ($T_f\in[1,1.3]$), 65 Sod samples ($T_f\in[1.5,2]$), and 84 sine-entropy samples ($T_f\in[1,1.2]$), each with polynomial degrees $p=\{1,2,3,4\}$. This yields 142, 85, and 91 discontinuity windows for the Lax, Sod, and sine-entropy problems, respectively, which are filtered as in Section~\ref{sec:hybrid-euler}. Additional details (e.g., filtering of the Euler entropy variable or different resolutions) can be found in Appendices~\ref{sec:entropy}--\ref{sec:multires}. Section~\ref{sec:results-datasets} first assesses the hybrid filter's performance across these datasets. Sections~\ref{sec:results-lax}--\ref{sec:results-shu-osher} then inspect final-time samples for each benchmark, focusing on discontinuity location and error reduction relative to the unfiltered $\mathcal{O}(1)$ error band. All error results are measured about the discontinuity windows in order to evaluate the performance of the hybrid filter for discontinuities. The performance of the hybrid filter in smooth regions has been demonstrated in prior SIAC theoretical work, e.g. \cite{Ryan2005}.

\subsection{Performance Assessment}\label{sec:results-datasets}

\begin{figure}[t]
    \centering
    \includegraphics[width=\linewidth]{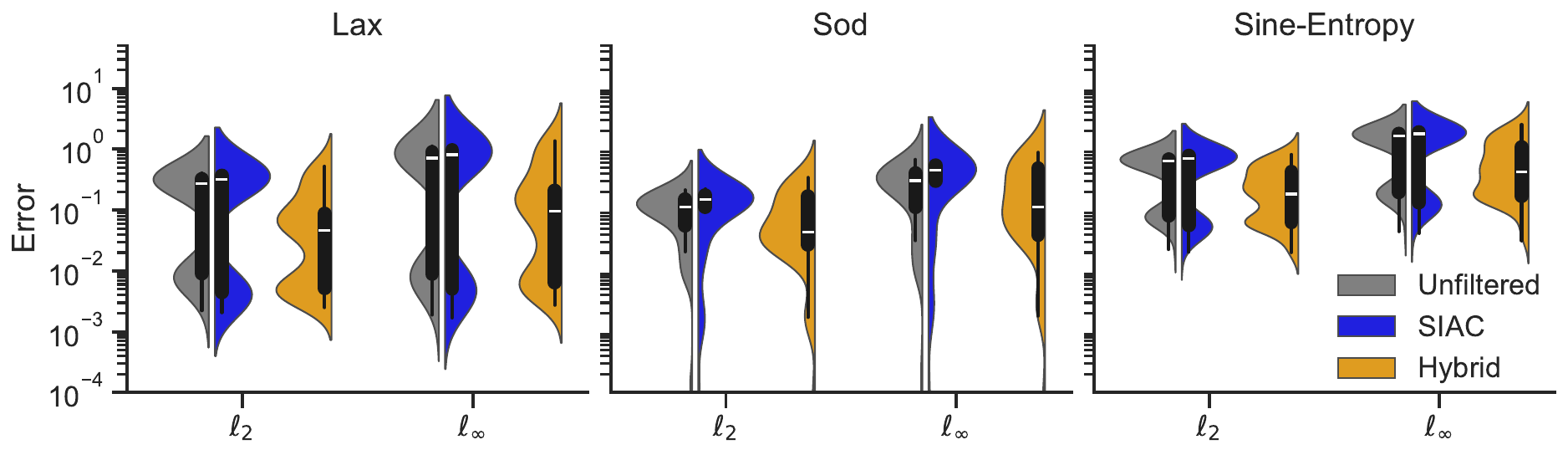}
    \caption{\small{\textbf{Euler velocity violinplots} for the grid $\ell_2$ and $\ell_\infty$ error in the \textit{(left to right) Lax, Sod, and sine-entropy (Shu-Osher)} datasets.}}
    \label{fig:dataset_velocity}
\end{figure}

\begin{figure}[t]
    \centering
    \includegraphics[width=\linewidth]{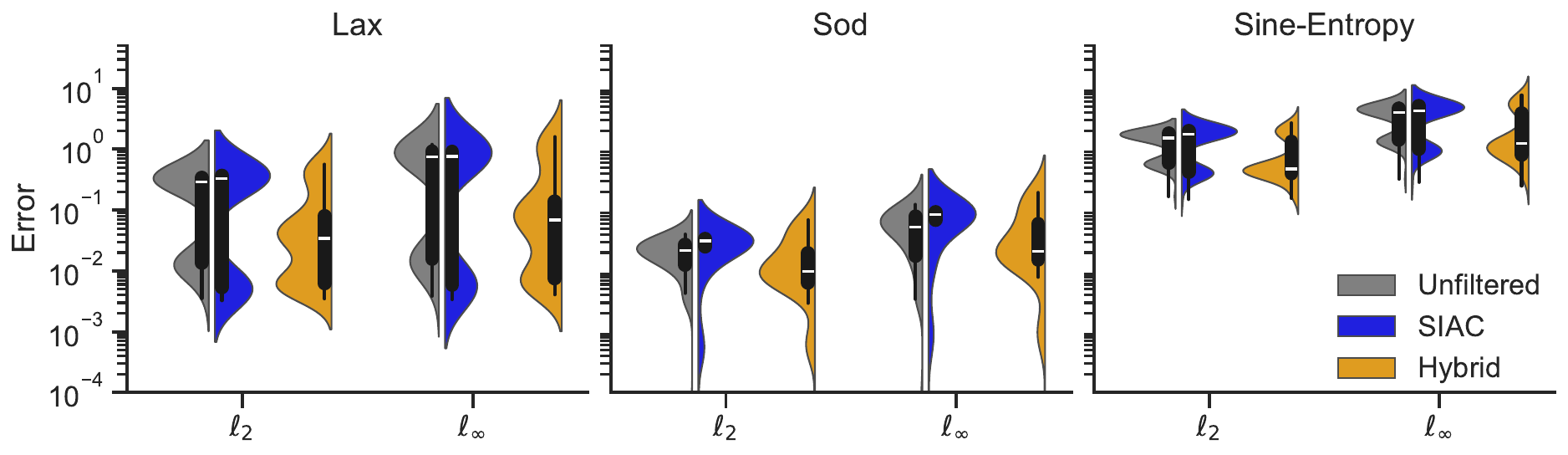}
    \caption{\small{\textbf{Euler pressure violinplots} for the grid $\ell_2$ and $\ell_\infty$ error in the \textit{(left to right) Lax, Sod, and sine-entropy (Shu-Osher)} datasets.}}
    \label{fig:dataset_pressure}
\end{figure}

\noindent Table~\ref{tab:error_dataset} shows Lax quartile errors for the unfiltered, SIAC-filtered, and hybrid-filtered approximations. Bold values mark hybrid errors below the unfiltered baseline, highlighting that the median error is consistently lower than the hybrid method. Figures~\ref{fig:dataset_density}--\ref{fig:dataset_pressure} compare error distributions for the Lax, Sod, and sine-entropy problems, with violinplots revealing lower error medians for the hybrid filter. Although Sod’s unfiltered errors are already moderate, the hybrid filter still offers a reduction. The sine-entropy problem is more complex, yet,  the hybrid filter, in particular, improves velocity and pressure about the discontinuity.

\subsection{Lax Samples at Final Time}\label{sec:results-lax}
\begin{figure}[t]
    \centering
    \includegraphics[width=0.7\linewidth]{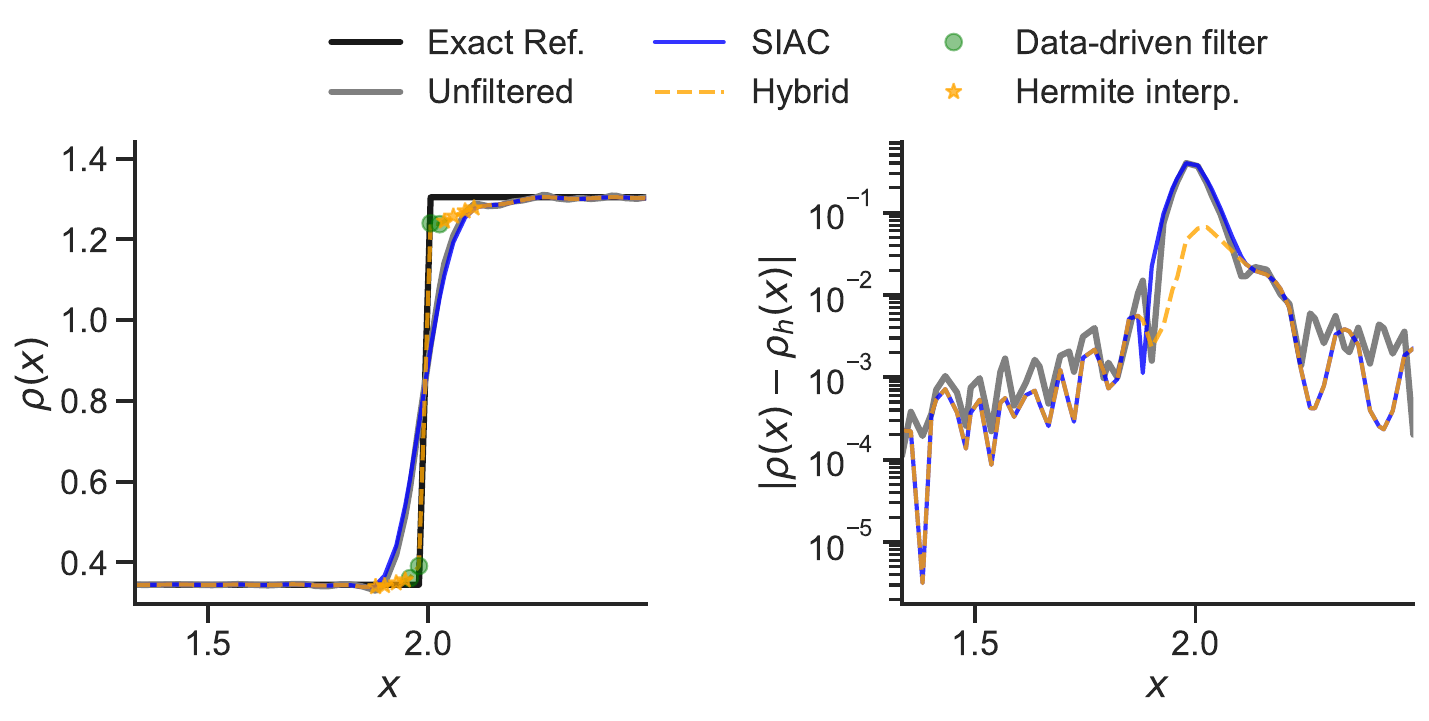}
    \caption{\small{\textbf{Lax density contact discontinuity} at $T_f=1.3$ (DG degree 2). \textit{Right:} unfiltered (gray), SIAC (blue), and hybrid (dashed yellow) solutions, with data-driven filter (green dots) and Hermite interpolation (yellow stars). \textit{Left:} pointwise errors.}}
    \label{fig:lax_contact_p2}
\end{figure}

\noindent Lax has two key density discontinuities near $x=2$ (contact) and $x=3.1$ (shock). Figure~\ref{fig:lax_contact_p2} shows that the hybrid filter pinpoints the contact, flattening the $\mathcal{O}(1)$ error band to $\mathcal{O}(10^{-1})$. Table~\ref{tab:lax_density} lists corresponding $\ell_2$/$\ell_\infty$ errors, with noticeable improvements for degrees $p>1$ at both discontinuities. The shock location is also well captured, as seen in Fig.~\ref{fig:lax_shock_p2}. Notably, degree 1 shock data benefit from up to an order-of-magnitude $\ell_2$ and $\ell_\infty$ error drop.

\subsection{Sod Samples at Final Time}\label{sec:results-sod}
\noindent For Sod, the contact discontinuity in density at $x\approx2$ does not appear in the velocity field, as it is constant there. Figure~\ref{fig:sod_contact_p2} shows that the hybrid filter preserves this constant velocity and reduces the error slightly (Table~\ref{tab:sod_velocity}). Around the shock ($x\approx-0.2$), Fig.~\ref{fig:sod_shock_p2} demonstrates the hybrid filter’s sharp location capture, cutting both $\ell_2$ and $\ell_\infty$ errors by nearly an order of magnitude for all $p$.

\subsection{Sine-Entropy (Shu-Osher) Samples at Final Time}\label{sec:results-shu-osher}
\noindent The sine-entropy (Shu-Osher) problem features a shock interacting with high-to-low frequency sinusoidal waves, which we demonstrate here for both the density and pressure profiles, Fig.~\ref{fig:so_shock_density_p2} and Fig.~\ref{fig:so_shock_pressure_p2}, respectively. The $\ell_2$ and $\ell_\infty$ errors included in Table~\ref{tab:sine-entropy} show that the SIAC filter performs better than the hybrid filter in reducing the $\ell_\infty$ error across all degrees and the unfiltered density data has the lowest $\ell_2$ error.  However, the hybrid-filtered data presented in Fig.~\ref{fig:so_shock_density_p2} is more accurate to the right of the shock, despite the relatively coarse mesh approximation (typically results are performed with $N=256$ elements). We expect the performance to be further improved with better resolution, as presented in the multi-resolution Lax filter results in Appendix~\ref{sec:multires}. The pressure profile in Fig.~\ref{fig:so_shock_pressure_p2} shows that the location of the discontinuity is resolved, and the error values in Table~\ref{tab:sine-entropy} indicate that the hybrid filter performs well in the $\ell_2$ and $\ell_\infty$ error for the degree 2 and 3 approximations while the SIAC filter improves the $\ell_\infty$ error for the other degrees.

\section{Discussion and Future Work} 
\label{sec:discussion}

\par \hypertarget{discuss}{In this work,} we introduce a hybrid post-processing filter that fuses an adapted SIAC kernel with data-driven and interpolation methods to reduce the $\mathcal{O}(1)$ error at discontinuities—a challenge that traditional filtering techniques have not overcome. It is only necessary to apply the filter at the final time. By training solely on top-hat data, the hybrid filter accurately locates sharp vertical slopes, precisely detecting shock and contact discontinuities in the Lax and Sod problems across density, velocity, and pressure. Although the shock in density for the Lax problem is only partially resolved, sub-cell feature extraction is achieved. In the sine-entropy (Shu-Osher) problem, our method notably improves the solution on one side of the shock despite the nonlinear dynamics and coarse simulation data. Overall, higher-order approximations benefit most, with a general narrowing of the $\mathcal{O}(1)$ error band near discontinuities.

\par Future work will refine the hybrid SIAC–data-driven filter for a broader problem class, including both low- and high-order approximations and enhanced precision for complex discontinuities. Moreover, while our current application is limited to post-processing at the final time, extending the filter to every time step could provide a more robust limiting strategy. Finally, we will explore optimization-based implicit neural network architectures~\cite{heaton2022explainable} to embed additional physical constraints into the filtering process.

\section*{Acknowledgments}
This work was partially funded by NSF award number DMS-2110745 and the U.S. Air Force Office of Scientific Research (AFOSR) Computational Mathematics Program (Program Manager: Dr. Fariba Fahroo) under Grant number FA9550-20-1-0166. This work is done in connection with Digital Futures and the Linné Flow Centre at KTH.

\bibliographystyle{siam}
\bibliography{refs}

\newpage 

\label{sec:Appendix}
\begin{appendices}

\begin{table}[t]
    \caption{\small{\textbf{Lax density} $\ell_2$ and $\ell_\infty$ errors for the \emph{contact} ($x\approx2$) and \emph{shock} ($x\approx3.1$) at $T_f=1.3$. Discontinuity windows are $[\min(S_i)-4,\ \max(S_i)+4]$.}}
    \centering
    \begin{tabular}{|c|c|c|c|c|c|}
    \hline 
         \multirow{2}{*}{Degree} & \multirow{2}{*}{Approx.} & \multicolumn{2}{c|}{Contact} & \multicolumn{2}{c|}{Shock} \\
         \cline{3-6} 
         & & $\ell_2$ & $\ell_\infty$ & $\ell_2$ & $\ell_\infty$ \\
          \hline 
       $p=1$ & Unfiltered & 2.62e-01 & 4.55e-01 & 1.50e-01 & 3.75e-01 \\
         & SIAC & 2.64e-01 & 4.55e-01 & 1.53e-01 & 3.89e-01 \\
         & Hybrid & 2.25e-01 & 7.30e-01 & 3.75e-02 & 6.78e-02\\
        \hline
        $p=2$ & Unfiltered & 1.93e-01 & 3.99e-01 & 1.44e-01 & 3.86e-01 \\
         & SIAC & 2.03e-01 & 3.97e-01 & 1.60e-01 & 4.15e-01 \\
         & Hybrid  & 3.99e-02 & 6.64e-02 & 1.06e-01 & 3.75e-01\\
         \hline 
         $p=3$ & Unfiltered & - & - & 1.50e-01 & 3.92e-01\\
         & SIAC & - & - & 1.60e-01 & 4.17e-01 \\
         & Hybrid & - & - & 1.20e-01 & 4.27e-01\\
         \hline 
         $p=4$ & Unfiltered & 1.84e-01 & 3.81e-01 & 1.53e-01 & 3.87e-01\\
         & SIAC & 1.94e-01 & 3.84e-01 & 1.63e-01 & 4.14e-01\\
         & Hybrid & 3.82e-02 & 5.88e-02 & 9.57e-02 & 3.37e-01 \\
         \hline 
    \end{tabular}
    \label{tab:lax_density}
\end{table}

\begin{figure}[t]
    \centering
    \includegraphics[width=0.7\linewidth]{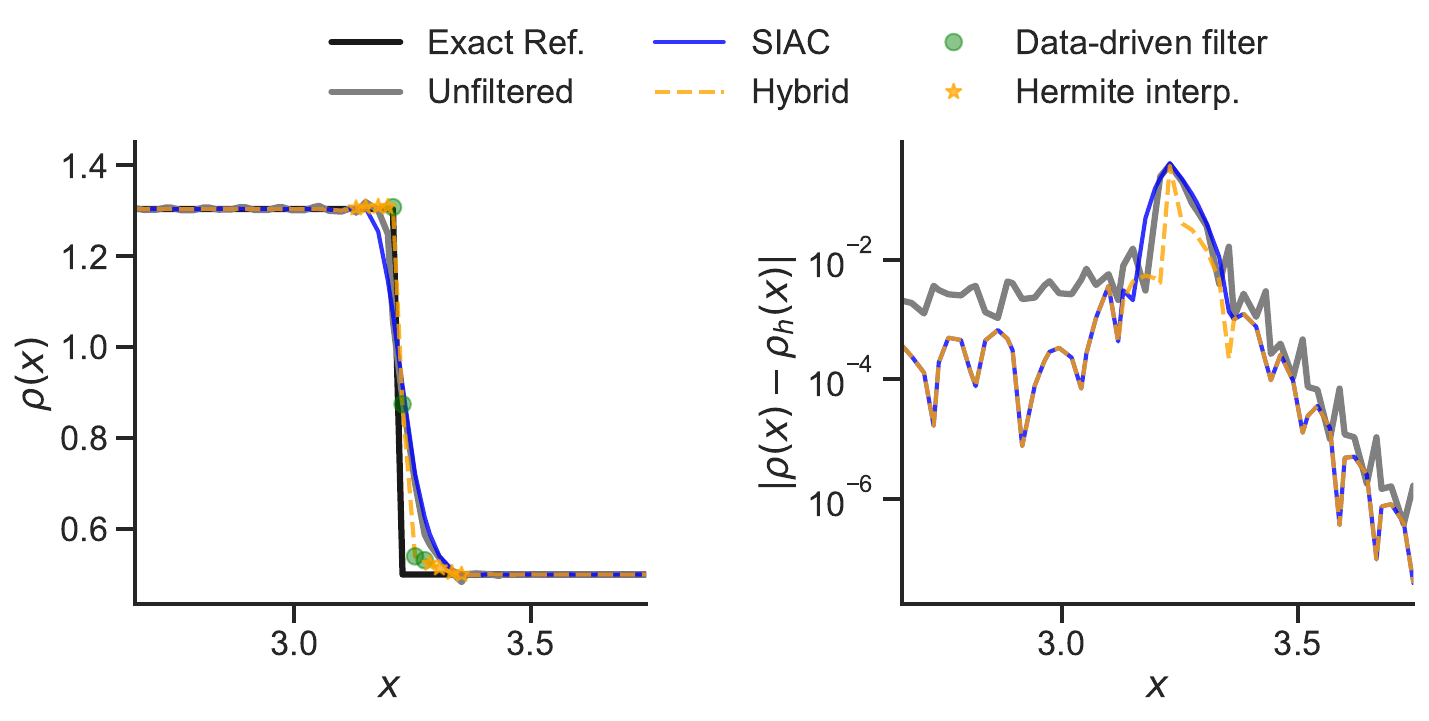}
    \caption{\small{\textbf{Lax density shock discontinuity} at $T_f=1.3$ (DG degree 2). \textit{Right:} unfiltered (gray), SIAC (blue), and hybrid (dashed yellow) solutions. \textit{Left:} pointwise errors.}}
    \label{fig:lax_shock_p2}
\end{figure}

\begin{figure}[t]
    \centering
    \includegraphics[width=0.7\linewidth]{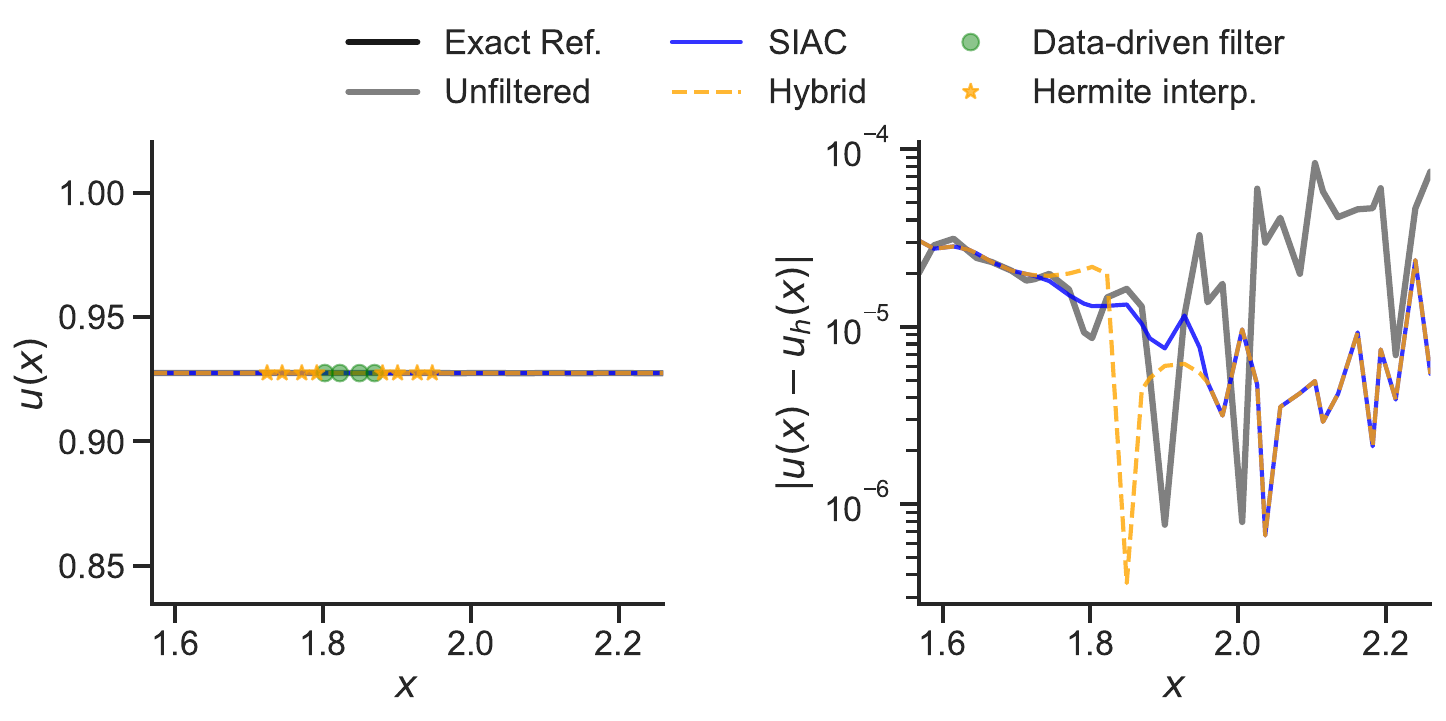}
    \caption{\small{\textbf{Sod velocity at the density contact discontinuity} at final time $T_f=2$ with DG degree 2 data: \textit{(right)} unfiltered (gray), SIAC moving average filter (blue), and hybrid filter (yellow dashed) approximations, including where the data-driven filter (green dots) and Hermite polynomial interpolation (yellow stars) are applied. The left plot shows the pointwise error for the respective filtered approximations. }}
    \label{fig:sod_contact_p2}
\end{figure}

\begin{table}[t]
    \caption{\small{\textbf{Sod velocity} grid $\ell_2$ and $\ell_\infty$ error for the \textit{density contact and shock discontinuities} at final time $T_f=2$ about discontinuity windows, $[\min(S_i) - 4, \ \max(S_i) + 4]$. \textit{Note: the Sod velocity profile has a jump discontinuity at the shock and not at the contact discontinuity for density; the error values at the contact are thus shown in italics.}}}
    \vspace{5pt}
    \centering
    \begin{tabular}{|c|c|c|c|c|c|}
    \hline 
        \multirow{2}{*}{Degree} & \multirow{2}{*}{Approx.}& \multicolumn{2}{c|}{Contact}& \multicolumn{2}{c|}{Shock} \\
         \cline{3-6} 
         & & $\ell_2$ & $\ell_\infty$& $\ell_2$ & $\ell_\infty$ \\
          \hline 
       $p=1$ & Unfiltered & \textit{5.73e-05} & \textit{6.88e-05} & 1.78e-01 & 5.37e-01\\
         & SIAC & \textit{5.49e-05} & \textit{6.43e-05} & 1.85e-01 & 4.80e-01\\
         & Hybrid & \textit{5.90e-05} & \textit{7.18e-05} & 3.66e-02 & 4.88e-02 \\ 
        \hline
        $p=2$ & Unfiltered &  \textit{5.78e-05} & \textit{8.37e-05} & 1.54e-01 & 4.04e-01\\
         & SIAC &  \textit{2.55e-05} & \textit{3.06e-05} & 1.87e-01 & 5.29e-01\\
         & Hybrid & \textit{2.63e-05} & \textit{3.06e-05} & 2.76e-02 & 5.84e-02\\
         \hline 
         $p=3$ & Unfiltered &   - & - & 1.50e-01 & 4.45e-01\\
         & SIAC &  - & - & 2.00e-01 & 5.60e-01\\
          & Hybrid &  - & - & 1.72e-02 & 4.69e-02\\
         \hline 
         $p=4$ & Unfiltered & \textit{8.58e-03} & \textit{1.02e-02} & 8.44e-02 & 2.67e-01\\
         & SIAC & \textit{2.21e-03} & \textit{3.03e-03} & 2.05e-01 & 5.71e-01\\
         & Hybrid & \textit{1.97e-03} & \textit{3.03e-03} & 1.69e-02 & 3.73e-02\\
         \hline 
    \end{tabular} 
    \label{tab:sod_velocity}
\end{table}

\begin{figure}
    \centering
    \includegraphics[width=0.7\linewidth]{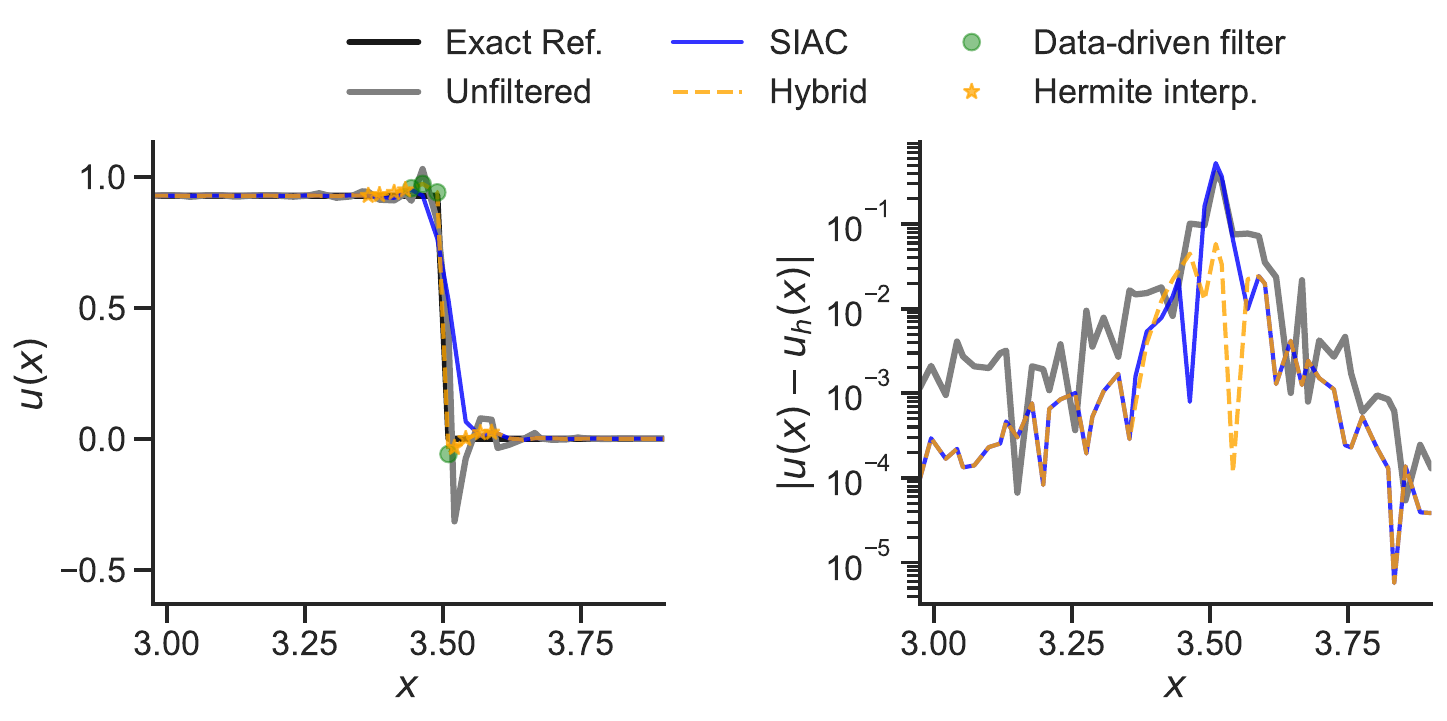}
    \caption{\small{\textbf{Sod velocity at the density shock discontinuity} at final time $T_f=2$ with DG degree 2 data: \textit{(right)} unfiltered (gray), SIAC moving average filter (blue), and hybrid filter (yellow dashed) approximations, including where the data-driven filter (green dots) and Hermite polynomial interpolation (yellow stars) are applied. The left plot shows the pointwise error for the respective filtered approximations.  }}
    \label{fig:sod_shock_p2}
\end{figure}
\section{Filtering the Entropy Variable}\label{sec:entropy}

The \hypertarget{entropy}{entropy variable} is another key Euler quantity that can be further examined and utilized for discontinuity detection \cite{Thea}. Similar to the primitive variables (velocity and pressure), filtering of the entropy variable can be achieved by first filtering the conservative variables and then evaluating the related quantities, as demonstrated in Section~\ref{sec:hybrid-euler}. The following entropy variable is considered, 
\[S = \frac{p}{\rho^\gamma}.\]
To examine the filtered entropy data, first filtering of density $\rho$, momentum $\rho u$, and energy $E$ is applied, then the filtered pressure $p=(\gamma-1)(E-\frac{1}{2}\rho u^2)$ is evaluated from the filtered conservative variables, and finally filtered entropy $S$ can be computed. The Lax entropy data is examined herein at the density contact and shock discontinuities.

The entropy data has an error profile of $\mathcal{O}(10)$ at the density contact discontinuity, as shown in the DG degree 2 data in Fig.~\ref{fig:lax_contact_entropy_p2} and in Table~\ref{tab:lax_entropy}. Despite the hybrid filter's capacity to handle smaller error profiles, the hybrid filtered entropy data shows the precise discontinuity location, as shown in Fig.~\ref{fig:lax_contact_entropy_p2}.

In contrast, the entropy data for the shock discontinuity is of $\mathcal{O}(10^{-1})$, 100$\times$
smaller than entropy at the density contact discontinuity. While the location of the discontinuity at the shock is extracted, there are more significant undershoots of the data, as shown in Fig.~\ref{fig:lax_shock_entropy_p2}, resulting in the larger $\ell_2$ and $\ell_\infty$ error, as reported in Table~\ref{tab:lax_entropy}. The SIAC filter performs better for most cases and \cite{PickloEdoh_Entropy} shows that it is promising for entropy correction.

\section{Multi-resolution Hybrid Filters}\label{sec:multires}

\hypertarget{MRAH}{Although the hybrid filter} was constructed for DG numerical data evaluated at $N=128$ elements in a domain $\Omega=[-5, 5]$, we test the performance of the hybrid filter at coarser ($N=64$) and finer ($N=256$) resolutions for the Lax problem. We expect the Sod problem to be similar, but note that the Lax problem poses more challenging discontinuities and so does the sine-entropy problem. Similar to Section~\ref{sec:results-lax}, we present here the filter results for the Lax density data at the contact and shock discontinuities. The discontinuity windows (Section~\ref{sec:disc_window}) for coarser data with $N=64$ elements were adapted to group troubled cells that are at most $n=2$ cells apart and with padding by $d=2$ cells on each side for a filter evaluation window of $[\min(S_i) - 2, \ \max(S_i) + 2]$. The discontinuity window parameters for finer data with $N=256$ elements were the same as those with $N=128$ elements ($n=d=4$). The hybrid filter construction is the same as defined in Section~\ref{sec:hybrid-app}.

For coarse data, the degree 2 filtered results for the contact discontinuity are presented in Fig.~\ref{fig:lax_N64_den_contact_p2} and the shock discontinuity are presented in Fig.~\ref{fig:lax_N64_den_shock_p2}. Both filtering outputs show an error reduction of the $\mathcal{O}(1)$ error band. The grid $\ell_2$ and $\ell_\infty$ errors (Table~\ref{tab:lax_N64_density}) are improved for degree 2. The hybrid filter improves the $\ell_2$ and $\ell_\infty$ error for the shock discontinuity for the degree 1-3 approximations but only improves the contact discontinuity in the $\ell_2$ error for degree 1 and in the $\ell_2$ and $\ell_\infty$ error for degree 2.

A resolution with $N=256$ elements is considered to examine the performance of the hybrid filter at finer resolutions. The degree 2 filtered data demonstrates precise contact (Fig.~\ref{fig:lax_N256_den_contact_p2}) and shock (Fig.~\ref{fig:lax_N256_den_shock_p2}) discontinuity location and feature extraction. Additionally, the $\ell_2$ and $\ell_\infty$ error values in Table~\ref{tab:lax_N256_density} demonstrate an error reduction for all degrees 1-4 where the discontinuities are detected.

\end{appendices}

\begin{figure}
    \centering
    \includegraphics[width=0.7\linewidth]{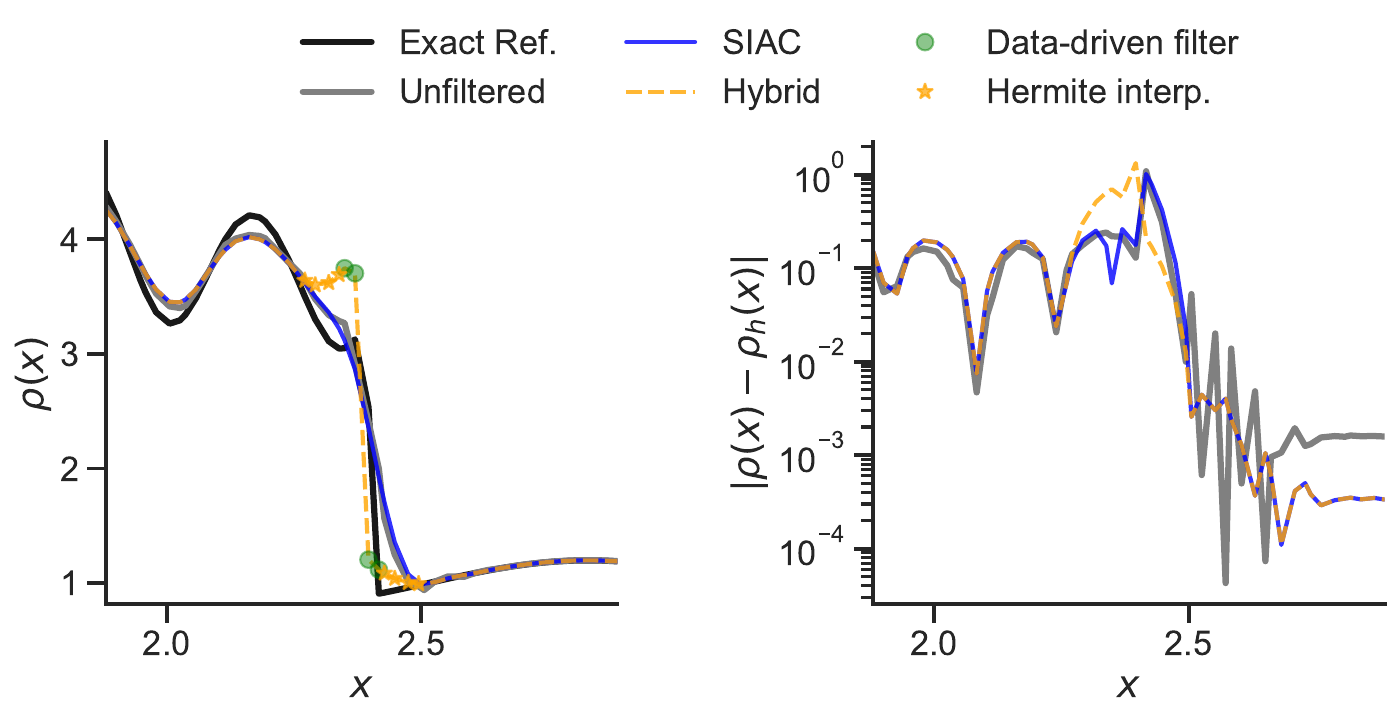}
    \caption{\small{\textbf{Sine-entropy (Shu-Osher) density at the density shock discontinuity} at final time $T_f=1.8$ with DG degree 2 data: \textit{(right)} unfiltered (gray), SIAC moving average filter (blue), and hybrid filter (yellow dashed) approximations, including where the data-driven filter (green dots) and Hermite polynomial interpolation (yellow stars) are applied. The left plot shows the pointwise error for the respective filtered approximations.  }}
    \label{fig:so_shock_density_p2}
\end{figure}

\begin{table}
    \caption{\small{\textbf{Sine-entropy (Shu-Osher) density and pressure} grid $\ell_2$ and $\ell_\infty$ error for the \textit{strong shock discontinuity} at final time $T_f=1.8$ about discontinuity window, $[\min(S_i) - 4, \ \max(S_i) + 4]$. }}
    \centering
    \begin{tabular}{|c|c|c|c|c|c|}
    \hline 
         \multirow{2}{*}{Degree} & \multirow{2}{*}{Approx.} & \multicolumn{2}{c|}{Density} & \multicolumn{2}{c|}{Pressure} \\
         \cline{3-6} 
         & & $\ell_2$ & $\ell_\infty$ & $\ell_2$ & $\ell_\infty$ \\
          \hline 
       $p=1$ & Unfiltered & 3.98e-01 & 9.75e-01 & 1.03e+00 & 3.21e+00\\
         & SIAC & 4.10e-01 & 8.40e-01 & 1.20e+00 & 3.10e+00\\
         & Hybrid & 5.81e-01 & 1.36e+00 & 1.24e+00 & 4.13e+00\\
        \hline
        $p=2$ & Unfiltered & 4.18e-01 & 1.09e+00 & 1.26e+00 & 3.77e+00\\
         & SIAC & 4.39e-01 & 1.02e+00 & 1.40e+00 & 3.66e+00\\
         & Hybrid  & 5.47e-01 & 1.33e+00 & 7.66e-01 & 2.35e+00 \\
         \hline 
         $p=3$ & Unfiltered & 4.08e-01 & 1.15e+00 & 1.22e+00 & 3.88e+00\\
         & SIAC & 4.16e-01 & 9.80e-01 & 1.36e+00 & 3.58e+00\\
         & Hybrid & 5.28e-01 & 1.24e+00 & 7.89e-01 & 2.66e+00 \\
         \hline 
         $p=4$ & Unfiltered & 3.80e-01 & 1.12e+00 & 1.12e+00 & 3.74e+00\\
         & SIAC & 4.00e-01 & 9.63e-01 & 1.31e+00 & 3.50e+00\\
         & Hybrid & 5.20e-01 & 1.23e+00 & 1.28e+00 & 4.51e+00\\
         \hline 
    \end{tabular}
    \label{tab:sine-entropy}
\end{table}

\begin{figure}
    \centering
    \includegraphics[width=0.7\linewidth]{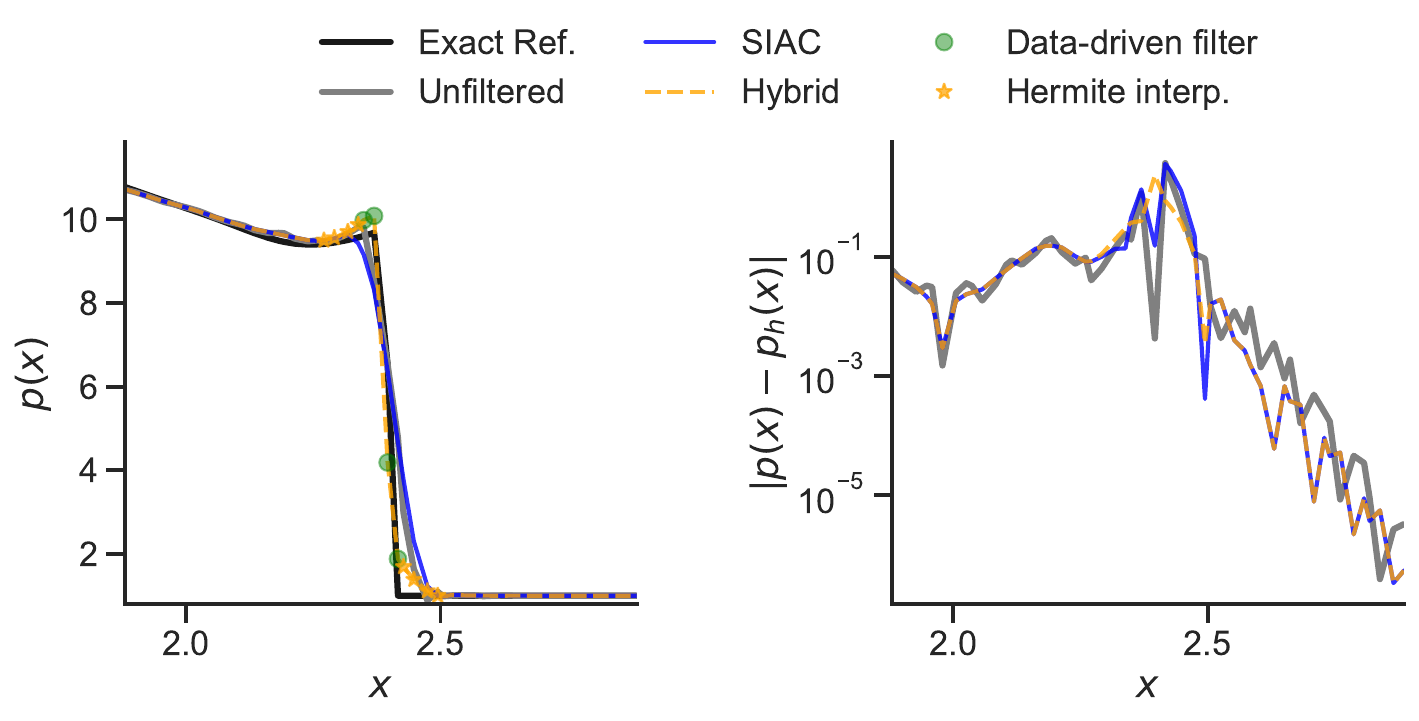}
    \caption{\small{\textbf{Sine-entropy (Shu-Osher) pressure at the density shock discontinuity}  at final time $T_f=1.8$ with DG degree 2 data: \textit{(right)} unfiltered (gray), SIAC moving average filter (blue), and hybrid filter (yellow dashed) approximations, including where the data-driven filter (green dots) and Hermite polynomial interpolation (yellow stars) are applied. The left plot shows the pointwise error for the respective filtered approximations.  }}
    \label{fig:so_shock_pressure_p2}
\end{figure}


\begin{figure}
    \centering
    \includegraphics[width=0.7\linewidth]{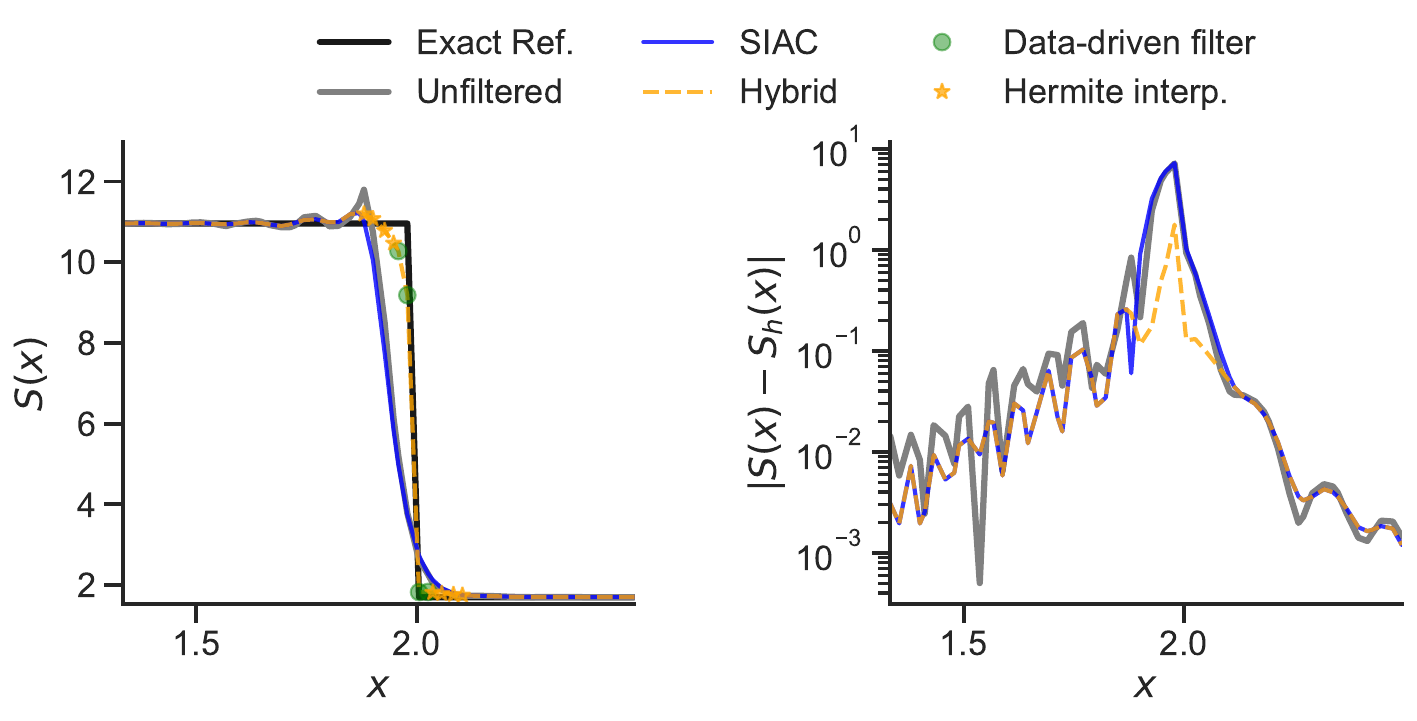}
    \caption{\small{\textbf{Lax entropy at the density contact discontinuity} at final time $T_f=1.3$ with DG degree 2 data: \textit{(right)} unfiltered (gray), SIAC moving average filter (blue), and hybrid filter (yellow dashed) approximations, including where the data-driven filter (green dots) and Hermite polynomial interpolation (yellow stars) are applied. The left plot shows the pointwise error for the respective filtered approximations.  }}
    \label{fig:lax_contact_entropy_p2}
\end{figure}

\begin{table}
    \caption{\small{\textbf{Lax $\boldsymbol{N=128}$ entropy} grid $\ell_2$ and $\ell_\infty$ error for the \textit{density contact and shock discontinuities} at final time $T_f=1.3$ about discontinuity windows, $[\min(S_i) - 4, \ \max(S_i) + 4]$. }}
    \centering
    \begin{tabular}{|c|c|c|c|c|c|}
    \hline 
         \multirow{2}{*}{Degree} & \multirow{2}{*}{Approx.}& \multicolumn{2}{c|}{Contact}& \multicolumn{2}{c|}{Shock} \\
         \cline{3-6} 
         & & $\ell_2$ & $\ell_\infty$& $\ell_2$ & $\ell_\infty$ \\
          \hline 
         $p=1$ & Unfiltered & 5.35e+00 & 7.52e+00 & 6.35e-02 & 1.81e-01 \\
         & SIAC & 5.00e+00 & 7.58e+00 & 9.96e-02 & 2.33e-01 \\
         & Hybrid & 3.99e+00 & 8.72e+00 & 6.85e-02 & 2.19e-01 \\
        \hline
        $p=2$ & Unfiltered & 3.02e+00 & 7.17e+00 & 7.30e-02 & 1.53e-01\\
         & SIAC & 3.15e+00 & 7.22e+00 & 8.81e-02 & 2.08e-01 \\
         & Hybrid & 5.71e-01 & 1.78e+00 & 2.18e-01 & 7.61e-01\\
         \hline 
         $p=3$ & Unfiltered & - & - & 7.50e-02 & 1.85e-01 \\
         & SIAC &  - & - & 9.39e-02 & 2.23e-01\\
         & Hybrid &  - & - & 2.24e-01 & 7.75e-01\\
         \hline 
         $p=4$ & Unfiltered & 2.83e+00 & 7.09e+00  & 7.48e-02 & 1.80e-01\\
         & SIAC & 2.98e+00 & 7.11e+00 & 9.12e-02 & 2.22e-01\\
         & Hybrid & 4.86e-01 & 1.68e+00 & 1.99e-01 & 6.88e-01 \\
         \hline 
    \end{tabular}
    \label{tab:lax_entropy}
\end{table}

\begin{figure}
    \centering
    \includegraphics[width=0.7\linewidth]{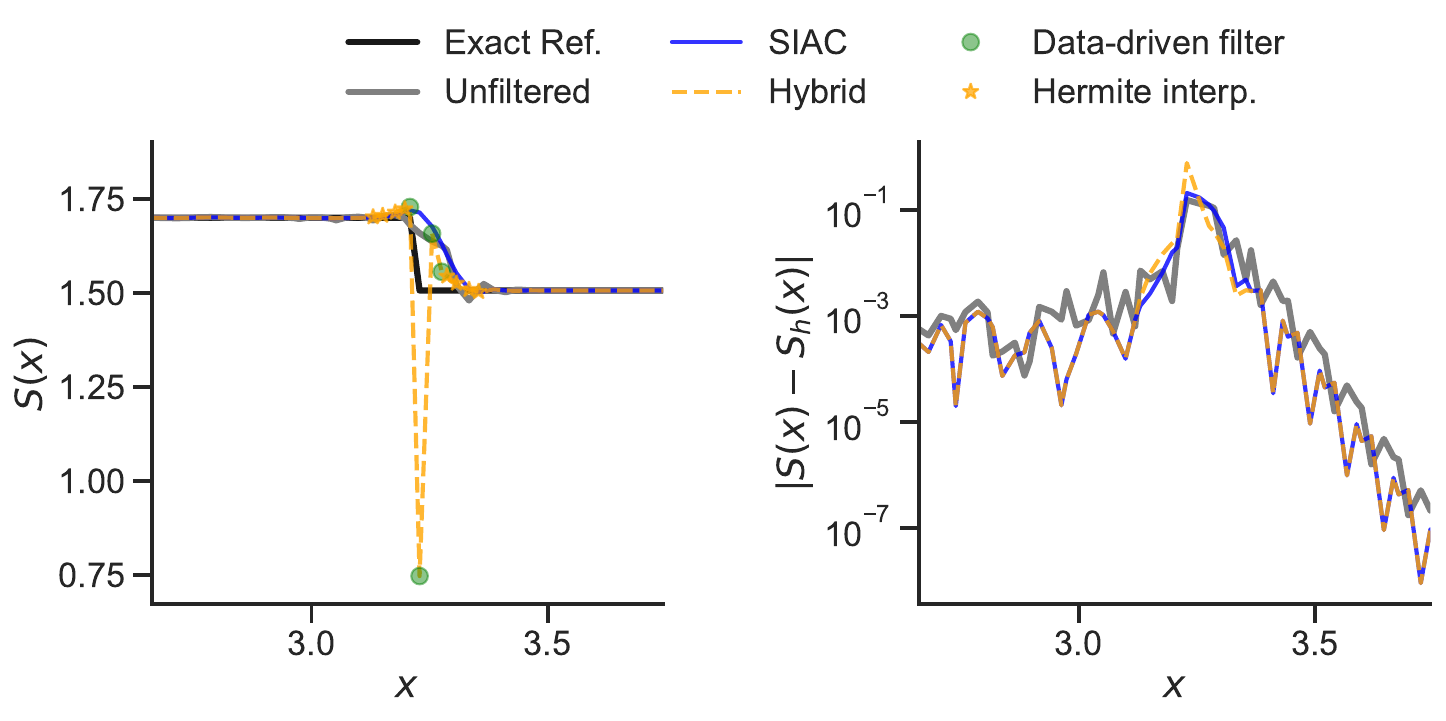}
    \caption{\small{\textbf{Lax entropy at the density shock discontinuity} at final time $T_f=1.3$ with DG degree 2 data: \textit{(right)} unfiltered (gray), SIAC moving average filter (blue), and hybrid filter (yellow dashed) approximations, including where the data-driven filter (green dots) and Hermite polynomial interpolation (yellow stars) are applied. The left plot shows the pointwise error for the respective filtered approximations. }}
    \label{fig:lax_shock_entropy_p2}
\end{figure}

\begin{figure}
    \centering
    \includegraphics[width=0.7\linewidth]{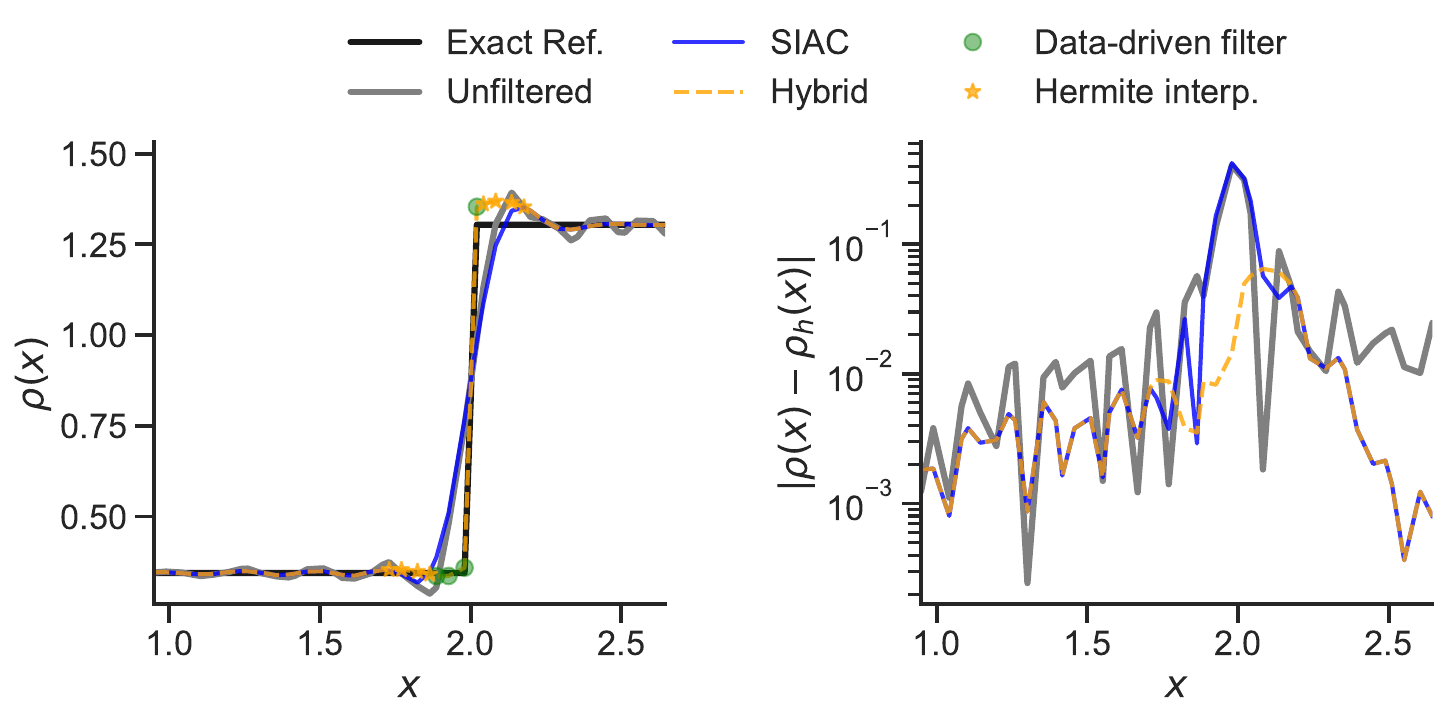}
    \caption{\small{\textbf{Coarser ($N=64$) Lax density contact discontinuity} at final time $T_f=1.3$ with DG degree 2 data: \textit{(right)} unfiltered (gray), SIAC moving average filter (blue), and hybrid filter (yellow dashed) approximations, including where the data-driven filter (green dots) and Hermite polynomial interpolation (yellow stars) are applied. The left plot shows the pointwise error for the respective filtered approximations.  }}
    \label{fig:lax_N64_den_contact_p2}
\end{figure}

\begin{table}
    \caption{\small{\textbf{Coarser ($N=64$) Lax density} grid $\ell_2$ and $\ell_\infty$ error for the \textit{density contact and shock discontinuities} at final time $T_f=1.3$ about discontinuity windows, $[\min(S_i) - 2, \ \max(S_i) + 2]$. }}
    \centering
    \begin{tabular}{|c|c|c|c|c|c|}
    \hline 
         \multirow{2}{*}{Degree} & \multirow{2}{*}{Approx.}& \multicolumn{2}{c|}{Contact}& \multicolumn{2}{c|}{Shock} \\
         \cline{3-6} 
         & & $\ell_2$ & $\ell_\infty$& $\ell_2$ & $\ell_\infty$ \\
          \hline 
          $p=1$ & Unfiltered & 3.48e-01 & 4.93e-01 & 1.91e-01 & 3.75e-01 \\
         & SIAC & 3.53e-01 & 4.94e-01 & 1.71e-01 & 3.47e-01 \\
         & Hybrid & 3.34e-01 & 7.93e-01 & 7.62e-02 & 9.91e-02\\
        \hline
        $p=2$ & Unfiltered & 2.30e-01 & 4.06e-01 & 1.85e-01 & 3.84e-01\\
         & SIAC & 2.40e-01 & 4.23e-01 & 1.55e-01 & 3.38e-01\\
         & Hybrid & 5.45e-02 & 6.46e-02 & 7.10e-02 & 1.12e-01\\
         \hline 
         $p=3$ & Unfiltered & - & - & 1.82e-01 & 3.80e-01  \\
         & SIAC &- & - &  1.71e-01 & 3.68e-01  \\
         & Hybrid & - & - &  1.13e-01 & 2.62e-01 \\
         \hline 
         $p=4$ & Unfiltered & 2.69e-01 & 4.72e-01 & 1.91e-01 & 4.61e-01\\
         & SIAC & 2.83e-01 & 4.50e-01  & 2.11e-01 & 4.39e-01 \\
         & Hybrid & 2.77e-01 & 6.81e-01 & 3.03e-01 & 7.62e-01\\
         \hline 
    \end{tabular}
    \label{tab:lax_N64_density}
\end{table}

\begin{figure}
    \centering
    \includegraphics[width=0.7\linewidth]{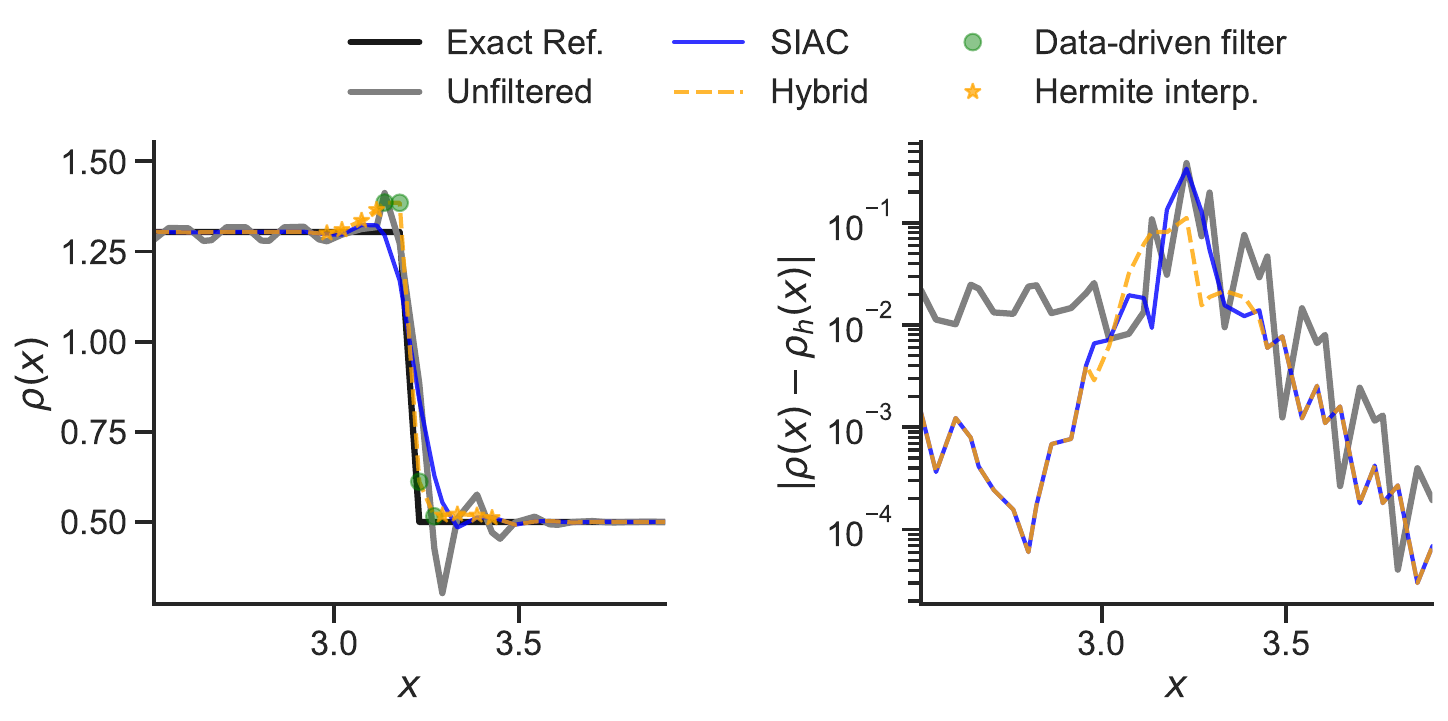}
    \caption{\small{\textbf{Coarser ($N=64$) Lax density shock discontinuity} at final time $T_f=1.3$ with DG degree 2 data: \textit{(right)} unfiltered (gray), SIAC moving average filter (blue), and hybrid filter (yellow dashed) approximations, including where the data-driven filter (green dots) and Hermite polynomial interpolation (yellow stars) are applied. The left plot shows the pointwise error for the respective filtered approximations.  }}
    \label{fig:lax_N64_den_shock_p2}
\end{figure}

\begin{figure}
    \centering
    \includegraphics[width=0.7\linewidth]{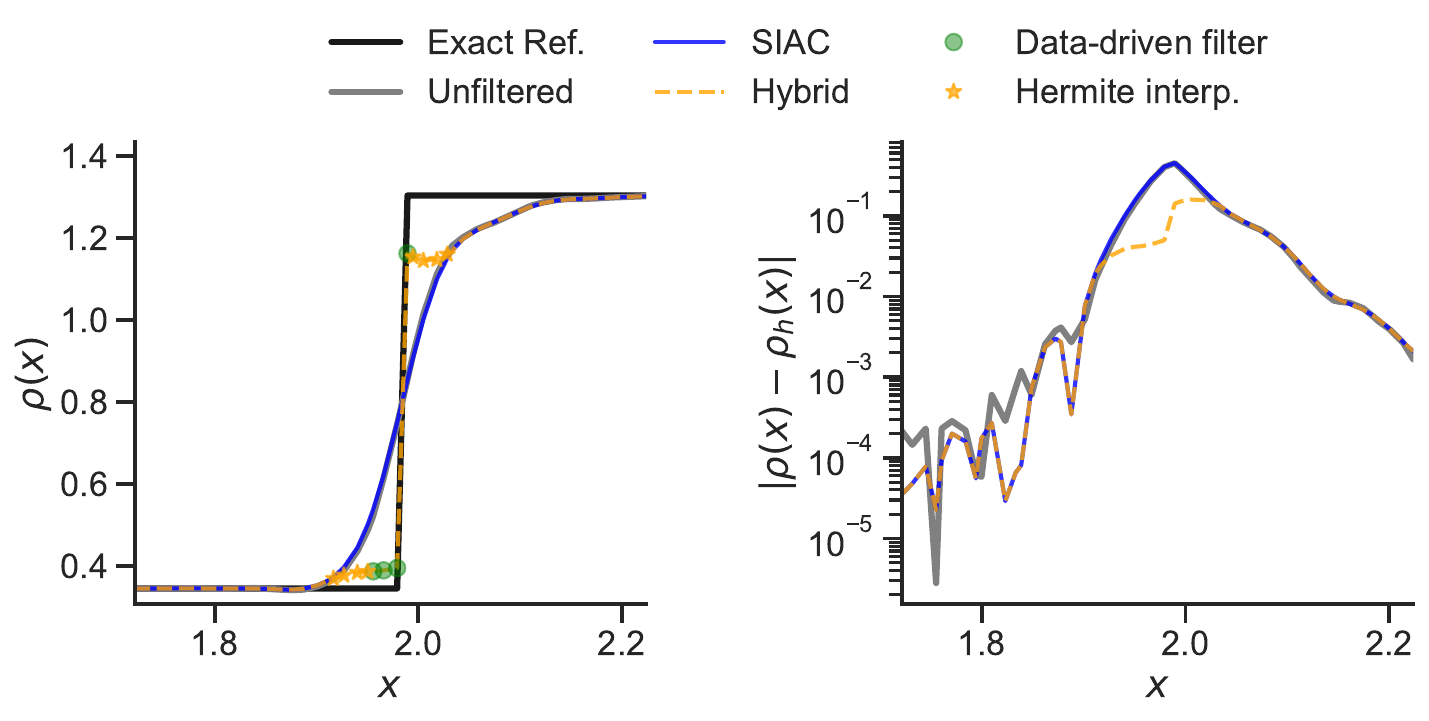}
    \caption{\small{\textbf{Finer ($N=256$) Lax density contact discontinuity} at final time $T_f=1.3$ with DG degree 2 data: \textit{(right)} unfiltered (gray), SIAC moving average filter (blue), and hybrid filter (yellow dashed) approximations, including where the data-driven filter (green dots) and Hermite polynomial interpolation (yellow stars) are applied. The left plot shows the pointwise error for the respective filtered approximations.   }}
    \label{fig:lax_N256_den_contact_p2}
\end{figure}

\begin{table}
    \caption{\small{\textbf{Finer ($N=256$) Lax density} grid $\ell_2$ and $\ell_\infty$ error for the \textit{density contact and shock discontinuities} at final time $T_f=1.3$ about discontinuity windows, $[\min(S_i) - 4, \ \max(S_i) + 4]$. }}
    \centering
    \begin{tabular}{|c|c|c|c|c|c|}
    \hline 
         \multirow{2}{*}{Degree} & \multirow{2}{*}{Approx.}& \multicolumn{2}{c|}{Contact}& \multicolumn{2}{c|}{Shock} \\
         \cline{3-6} 
         & & $\ell_2$ & $\ell_\infty$& $\ell_2$ & $\ell_\infty$ \\
          \hline 
          $p=1$ & Unfiltered & 2.11e-01 & 4.82e-01 & 1.07e-01 & 4.18e-01\\
         & SIAC & 2.14e-01 & 4.83e-01 & 1.11e-01 & 3.98e-01\\
         & Hybrid & 1.33e-01 & 3.02e-01 & 3.83e-02 & 1.77e-01 \\
        \hline
        $p=2$ & Unfiltered & 1.81e-01 & 4.49e-01 & 1.09e-01 & 4.14e-01 \\
         & SIAC & 1.86e-01 & 4.50e-01 & 1.17e-01 & 4.05e-01 \\
         & Hybrid & 8.33e-02 & 1.59e-01 & 4.45e-02 & 2.09e-01 \\
         \hline 
         $p=3$ & Unfiltered & - & - & 1.07e-01 & 4.12e-01\\
         & SIAC &  - & - & 1.15e-01 & 4.07e-01\\
         & Hybrid &  - & - & 4.98e-02 & 2.37e-01 \\
         \hline 
         $p=4$ & Unfiltered & 1.78e-01 & 4.54e-01 & 1.10e-01 & 4.16e-01 \\
         & SIAC & 1.84e-01 & 4.60e-01 & 1.16e-01 & 4.08e-01\\
         & Hybrid & 8.42e-02 & 1.65e-01 & 4.92e-02 & 2.33e-01 \\
         \hline 
    \end{tabular}
    \label{tab:lax_N256_density}
\end{table}

\begin{figure}
    \centering
    \includegraphics[width=0.7\linewidth]{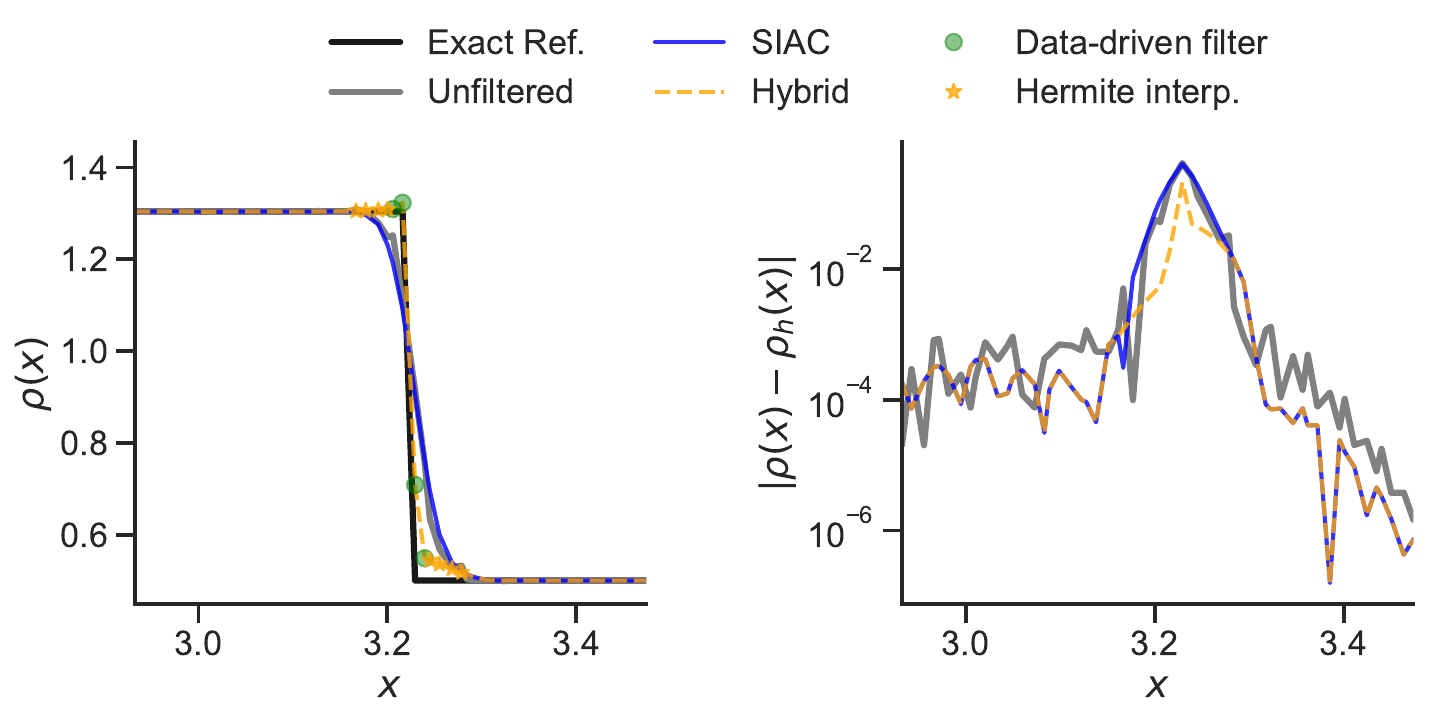}
    \caption{\small{\textbf{Finer ($N=256$) Lax density shock discontinuity} for DG degree 2 data: \textit{(right)} unfiltered (gray), SIAC moving average filter (blue), and hybrid filter (yellow dashed) approximations, including where the data-driven filter (green dots) and Hermite polynomial interpolation (yellow stars) are applied. The left plot shows the pointwise error for the respective filtered approximations. }}
    \label{fig:lax_N256_den_shock_p2}
\end{figure}

\end{document}